\documentclass[11pt]{amsart} \textwidth=14.5cm \oddsidemargin=1cm
\usepackage{youngtab}

\usepackage[latin1]{inputenc}
\usepackage{amsxtra}
\usepackage{amsmath}
\usepackage{amscd}
\usepackage{amssymb}
\usepackage{amsfonts}
\usepackage[all]{xy}
\usepackage{mathrsfs}
\usepackage{amsthm}
\usepackage{color}
\usepackage{upgreek}
\usepackage{verbatim}
\usepackage{enumerate}
\usepackage{latexsym}
\usepackage{graphicx}
\usepackage{tikz}
\usepackage{todonotes}
\usepackage{setspace}
\usepackage{hyperref}
\usepackage{stmaryrd}
\usepackage{nicefrac}
\usepackage{euscript}
\usetikzlibrary{decorations.pathmorphing}

 \definecolor{darkgreen}{HTML}{336633}
 \definecolor{darkred}{HTML}{993333}
\makeatletter

\newcommand{\arxiv}[1]{\href{http://arxiv.org/abs/#1}{\tt
    arXiv:\nolinkurl{#1}}}

\theoremstyle{plain}
\newtheorem{thm}{Theorem}[section]
\newtheorem*{thmA}{Theorem B}
\newtheorem*{thmB}{Theorem D}
\newtheorem*{thmC}{Theorem A}
\newtheorem*{thm99}{Theorem C}
\newtheorem{lem}[thm]{Lemma}
\newtheorem{prop}[thm]{Proposition}

\newtheorem{cor}[thm]{Corollary}
\newtheorem{df-prop}[thm]{Definition-Proposition}

\theoremstyle{definition}

\theoremstyle{remark}
\newtheorem{rem}[thm]{Remark}

\xyoption{all}
\def\Res{\operatorname{Res}\nolimits}
\def\Ind{\operatorname{Ind}\nolimits}

\def\ep{\epsilon}
\def\gl{\mathfrak{gl}}

\def\la{\lambda}

\def\pn{\mf{pe} (n)}
\def\ov{\overline}

\newcommand{\mc}{\mathcal}
\newcommand{\mf}{\mathfrak}
\newcommand{\C}{\mathbb C}
\newcommand{\oo}{{\ov 0}}
\newcommand{\oa}{{\bar 0}}
\newcommand{\ob}{{\bar 1}}
\newcommand{\vare}{\epsilon} 

\newcommand{\fg}{\mathfrak{g}}

\newcommand{\fb}{\mathfrak{b}}

\newcommand{\fh}{\mathfrak{h}}

\newcommand{\mZ}{\mathbb{Z}}
\newcommand{\cO}{\mathcal{O}}

\newcommand{\mC}{\mathbb{C}}

\newcommand{\h}{\mathfrak{h}}

\newcommand{\ch}{\mathrm{ch}}

\newcommand{\Coind}{{\rm Coind}}

\newcommand{\g}{\mathfrak{g}}

\newcommand{\fp}{\mathfrak{p}}

\newcommand{\nb}{{\nabla}}
\newcommand{\Z}{{\mathbb Z}}

\title[Parabolic category $\mc O^\fp$ for  periplectic Lie superalgebras $\pn$]{Parabolic category $\mc O^\fp$ for  periplectic Lie superalgebras $\pn$}

\author[Chen]{Chih-Whi Chen}
\address{Department of Mathematics, National Central University, Zhongli District, Taoyuan City, Taiwan 32001} 
\email{cwchen@math.ncu.edu.tw} 

\author[Peng]{Yung-Ning Peng}
\address{Department of Mathematics, National Central University, Zhongli District, Taoyuan City, Taiwan 32001}
\email{ynp@math.ncu.edu.tw}

\begin{document}
	\begin{abstract}We provide a linkage principle in  an arbitrary  parabolic category $\mc O^\fp$  for the periplectic Lie superalgebras $\pn$. As an application, we classify indecomposable  blocks in   $\mc O^\fp$. 
		We classify indecomposable tilting modules in $\mc O^\fp$ whose characters are controlled by the Kazhdan-Lusztig polynomials of type $\bf A$ Lie algebras. We establish the complete list of characters of indecomposable tilting modules in $\mc O^\fp$ for $\mf{pe}(3).$ 
	\end{abstract}

\numberwithin{equation}{section}

\maketitle


\section{Introduction}  \label{Sect::Introduction}
\subsection{} 

Character problems are central in representation theory. Since Kac's pioneering work in \cite{Kac77},  the irreducible characters of many finite-dimensional classical Lie superalgebras have been studied in connection with various areas of classical Lie theory (see, e.g., \cite{Br1, Br3, BW, CLW1, CLW2, EhS, GS, Se96})  with the exception of the periplectic Lie superalgebras $\pn$. 

The periplectic Lie superalgebra $\pn =\pn_\oa \oplus \pn_\ob$ is a superanalogue of the orthogonal or symplectic Lie algebra. There is a natural  $\Z$-grading   $\pn =\pn_{-1}\oplus \pn_0 \oplus \pn_1$ which is compatible with $\Z_2$-gradation such that $$\pn_0=\pn_\oa \cong \gl(n),$$ and
$$\pn_{-1}\cong \Lambda^2(\C^{n\ast})~\text{ and }\pn_{1} \cong S^2(\C^{n}),$$ as $\pn_\oa$-modules. 

Recently, a breakthrough in the study of the category $\mc F_n$ consisting of finite-dimensional modules for $\pn$ was achieved using the fake Casimir element (cf. \cite{B+9, Co}). Subsequently, the problem of irreducible character formulae in $\mc F_n$ has been solved in \cite{B+9}  by determining multiplicities of standard and costandard modules in indecomposable projective.
This has promoted a resurgence of interest in the representation theory of $\pn$ (see also \cite{Ch} for partial results). 

\subsection{}  
The description of  irreducible modules of blocks in   $\mc F_n$   has been independently obtained in \cite[Theorem 8.3.1]{Co} and \cite[Theorem 9.1.2]{B+9} (see also \cite{Ch} for partial results). In particular, the number of (indecomposable) blocks in $\mc F_n$ is $n+1$.  Later on, it is proved \cite[Theorem~5.4]{CC} that the  BGG category $\mc O_\Z$ of modules of integral weights for $\pn$ owns the same number of blocks.   

While the category $\mc F_n$ and the BGG category $\mc O$ for $\pn$ have been  extensively studied in recent years (see, e.g., \cite{B+9}, \cite{B+92}, \cite{BK}, \cite{CC}, \cite{Ch}, \cite{Co},    \cite{ES1}, \cite{ES2}, \cite{ES3}, \cite{Go}, \cite{HIR},  \cite{IRS}, \cite{Mo},  \cite{Se02}), however, the study about parabolic version of BGG category for $\pn$ remains to be unavailable in the literature.

The parabolic subalgebras of $\pn$ are   classified in terms of combinatorics of  bi-partitions in \cite[Section 5.3]{CCC}. It is natural to study the representation theory of the corresponding parabolic BGG categories.  The goal of the present paper is to study representations in an arbitrary parabolic BGG category for $\pn$. In this sequel to \cite{CC} and \cite{CCC}, we consider block decomposition and character problem in these categories.

  \subsection{} Throughout this article, we let $\mf g$ denote the periplectic Lie superalgebra $\pn$. Also, we fix a standard Cartan subalgebra  $\mf h$ and a standard Borel subalgebra $\mf b_\oa$  for $\g_\oa$ as defined in Section \ref{Intropen}.
This gives rise to a  distinguished  Borel subalgebra $\mf b:= \mf b_\oa \oplus \g_1$ of $\g$. The BGG category $\mc O \equiv \mc O(\g,\mf b)$  for $\g$ consists of $\mf h$-semisimple, finitely-generated $\g$-modules on which $\mf b$ acts locally finitely.

 Let $\{\vare_1, \vare_2, \ldots, \vare_n\}\subseteq \mf h^*$ denote the dual basis for the standard basis of $\mf h$ (see Section \ref{Intropen}).   The sets of even roots, positive even roots and odd roots of $\g$ are respectively denoted by $\Phi_\oa, \Phi^+_\oa$ and $\Phi_\ob$:
  	\begin{align}
  	&\Phi_\oa = \{\epsilon_i-\epsilon_j\,|\, 1\le i\not=j\le n\},~\Phi_\oa^+ = \{\epsilon_i-\epsilon_j\,|\, 1\le i<j\le n\},\label{Roots1}\\
  	&\Phi_\ob = \{- \epsilon_i-\epsilon_j\,|\, 1\le i< j\le n\} \cup  \{ \epsilon_i+\epsilon_j\,|\, 1\le i\le j\le n\}. \label{Roots2}
  	\end{align} In particular, the set $$\{\alpha_i:= \vare_i -\vare_{i+1}|~1\leq i\leq n-1\}$$ forms the simple system of $\Phi_\oa$. 

The Weyl group $W$ of $\g$ is defined to be the Weyl group of $\mf g_\oa$. We will identify $W$ with $\mf S_n$, the group generated by simple reflections $\{ s_{\alpha_i} \,|\, 1\leq i\leq n-1\}$.
We normalize the non-degenerate $\mf S_n$-invariant bilinear form $\langle\cdot, \cdot\rangle: \mf h^*\times \mf h^* \rightarrow \C$ by $\langle \vare_i, \vare_j \rangle  =\delta_{ij}$ for all $1\leq i,j\leq n$. It enables one to define usual notions in highest weight theory, e.g., integral, non-integral, dominant, anti-dominant weights.   Furthermore, for each positive even root $\beta = \vare_i -\vare_j \in \Phi_\oa^+$, we define the associated  conjugate $\ov  \beta$ to be the odd root $\ov \beta:=\vare_i+\vare_j \in \Phi_\ob$.
  
     Throughout, we follow the notions of parabolic decomposition given in \cite{CCC}  (see also  \cite{Ma}). We consider a (reduced) parabolic subalgebra $\fp$ satisfying $\mf b\subseteq \fp\subseteq \mf g$ with a purely even Levi subalgebra $\mf l\subseteq \fp$ as defined  in \cite[Section 1]{CCC}.  The even part $\mf p_\oa$ of $\mf p$ is a parabolic subalgebra of $\mf g_\oa$ arising from a parabolic decomposition of $\mf g_\oa$. 
The corresponding {\em parabolic category $\mc O^\fp :=\mc O(\mf g,\mf p)$} is a highest weight category with standard objects $\Delta^\fp_\la$, i.e. parabolic Verma modules, costandard objects  $\nb^\fp_\la$ and irreducible objects $L_\la$ indexed by  $\fp$-dominant (cf. \cite[Section 3]{CCC}) weights $\la\in \Sigma_\fp^+$; see Section 3.1 for the precise definition of the set $\Sigma_\fp^+$. In particular,  $L_\la$ is the irreducible module of highest weight  $\la -\rho$ with respect to the Borel subalgebra $\mf b$, where $\rho$ denotes the Weyl vector of $\g_\oa$.

    The Weyl group $W$ acts naturally on $\mf h^\ast$. Following \cite[Section 5.1]{Hu08}, we define the notion of strongly linked weights as follows. Let $\la,\mu\in \h^\ast$. We say $\mu$ is {\em strongly linked} to $\la$, denoted by $\mu\uparrow \la$, if $\la =w\mu$ for some $w\in W$ such that  either $w=e$ or there exist positive roots $\beta_1, \beta_2,\ldots , \beta_k\in \Phi_\oa^+$ such that  $w=s_{\beta_1}\cdots s_{\beta_k}$, $\langle\mu,\beta_k\rangle \in \Z_{>0}$ and $\langle s_{\beta_\ell}\cdots s_{\beta_k} \mu,~ \beta_{\ell-1}\rangle \in \mathbb Z_{> 0}$, for any $2\leq \ell\leq k$.

Let $  \Phi^+(\mf l)$ denote the set of positive roots of $\mf l$. The following is the first main result in this article, which enables one to construct linkages in $\mc O^\fp$. 

\begin{thmC} \label{thm::thmC}Let $\la \in \Sigma^+_\fp.$ We have the following sufficient conditions for positive multiplicities in parabolic  (dual) Verma modules $\Delta_\la^\fp$ and $\nb^\fp_{\la}$: 
	\begin{itemize}
		\item[(1)] Let $1\leq q\leq n$ be given such that $\la -2\vare_q\in \Sigma_\fp^+$.  If $\la-2\vare_q$ is not strongly linked to $\la -\ov \alpha_i$ for any $1\leq i\leq n-1$ with $\langle\la, \alpha_i\rangle =1$, then   
		$[\nb^\fp_{\la}:L_{\la-2\vare_q}]>0.$
		\item[(2)] Let  $\alpha_i\in \Phi^+(\mf l)$ be given such that $\la -\ov \alpha_i\in \Sigma_\fp^+$ and  $\langle\la, \alpha_i\rangle =1$. If $\la -\ov \alpha_i$ is not strongly linked to $\la-2\vare_q$ for any $1\leq q\leq n$, then   $[\Delta^\fp_{\la}:L_{\la-\ov \alpha_i}]>0.$
	\end{itemize}	
\end{thmC}

It was shown in  \cite{Br4}   (see also  \cite[Section 5.1]{CoM2} and \cite[Section 4]{St}) that each (integral) block in an arbitrary  parabolic category for the general linear Lie algebra $\gl(n)$ remains indecomposable whenever it is non-zero. As an application of Theorem A, we establish in Section \ref{Sect::parabolic}  a similar phenomenon for the periplectic Lie superalgebra $\pn$. The second  main result in this article is the following.
\begin{thmA} \label{thm::thmA} Each block in $\mc O$ remains indecomposable   when restricted to an arbitrary parabolic subcategory $\mc O^\fp$.	In particular, the number of blocks in the parabolic category $\mc O_\Z^\fp$ of modules of integral weights for   $\pn$ is $n+1$.
\end{thmA}

\subsection{} For an $\mf h$-semisimple module $M$ and a weight $\mu \in \h^\ast$,  let  $M^\mu:= \{m \in M~| h  m =\mu(h) m, \text{ for all } h\in \h\}$
 denote its $\mu$-weight space. Suppose that $M^\mu$ are finite-dimensional for all $\mu \in \h^\ast$. The {\em character} of $M$  is  defined to be the formal sum
\begin{align*}
\text{ch}M : = \sum_{\mu\in \h^\ast}\text{dim}M^\mu e^{\mu},\end{align*}
where $e$ is a formal indeterminate. 

The solution to the irreducible character problem in the BGG category for finite-dimensional semisimple Lie algebras has been given in Kazhdan-Lusztig polynomials, see \cite{BB, BK, KL}. Subsequently, the parabolic analogue of Kazhdan-Lusztig polynomials and their representation theoretic interpretations were   obtained,  see, e.g., \cite{BGS}, \cite{De} and \cite{So}.

For  $\la \in \Sigma_\fp^+$, we denote by $M^{\fp_\oa}_{\la}$   the parabolic Verma $\mf g_\oa$-module of highest weight $\la-\rho$ (see \eqref{Eq::paraV}). Also, let $L^\oa_\la$ denote the unique simple quotient of the $\mf g_\oa$-module  $M^{\fp_\oa}_{\la}$. In particular, the multiplicity $[M^{\fp_\oa}_\la:L_\la^\oa]$ can be computed by (parabolic) Kazhdan-Lusztig polynomials of type A Lie algebras, see, e.g., \cite{So}.

\subsection{} Recently, a version of Ringel duality for an arbitrary classical Lie superalgebra has been established in \cite[Theorem 3.7]{CCC} (see also \cite{BS18}), including $\pn$.    The character of an indecomposable tilting module $T^\fp_\la$ of highest weight $\la -\rho$ can be expressed in terms of  those of irreducible  modules 
\begin{align}
&\ch T^\fp_\la = \sum_{\mu \in \Sigma_\fp^+}  [\nabla^\fp_{-w_0^\fp\mu}: L_{-w_0^\fp\la}] \ch \Delta^\fp_\mu,
\end{align} where $w_0^\fp$ is the longest element in the parabolic subgroup of $W$ associated to the Levi subalgebra $\mf l.$

 The irreducible and tilting characters of finite-dimensional modules for $\pn$ were computed in \cite{B+9}, where an algorithm was given (purely combinatorially in terms of certain diagrams). 
However, it is highly non-trivial to relate the algorithm with Lusztig's canonical bases and classical Lie theory.
It is therefore surprising that there exists a large class of tilting modules whose characters can be obtained by the Kazhdan-Lusztig polynomials  of type A Lie algebras, showing a connection between the representation theory of $\pn$ and classical Lie theory.
  
To describe this class of modules, we recall that a weight $\la \in \h^\ast$ is said to be {\em typical} \cite[Section 5]{Se02}  if $$\prod_{\beta \in \Phi_\oa}( \langle \la, \beta\rangle   -1 ) \neq 0.$$
 In the present paper, we propose the following notion. A weight $\la\in \h^\ast$ is said to be {\em $\fp$-weakly-typical} if 
 \begin{align*}
 &\prod_{\beta \in \Phi^+(\mf l)} (\langle \la, \beta\rangle-1) \cdot \prod_{\gamma \in \Phi_\oa^+ \backslash \Phi^+(\mf l)}(\langle \la, \gamma\rangle+1)\neq 0.
 \end{align*}

  Our third main theorem shows that the character formula of tilting modules of $\mf p$-weakly-typical highest weights are given by (parabolic)  Kazhdan-Lusztig polynomials of type A Lie algebras:
 \begin{thm99}  \label{thm::thmD} If $\la$ $\in\Sigma^+_\fp$ is $\fp$-weakly-typical, then 
    \begin{align}
    &\ch T^\fp_\la =\sum_{\mu \in \h^\ast} [M^{\fp_\oa}_{-w_0^\fp\mu}:L^{\oa}_{-w_0^\fp\la}] \ch \nb^\fp_\mu.
    \end{align}
 \end{thm99}
 
 In particular, if $\mf p = \mf b$  and  $\la\in \h^\ast$ is a dominant integral weight, then $\la$ is $\mf b$-weakly-typical. Consequently, Theorem~C deduces the tilting character formula 
$$\ch T^{\mf b}_\la = \sum_{\mu\in W\la} \ch \nb^{\mf b}_\mu$$ 
in the full category $\mc O$.

  \subsection{} 
    The characters of projective covers in the full category $\mc O$  for $\mf{pe}(2)$ has been computed in \cite[Subsection 6.2]{CC}.  {Later on, the tilting characters are obtained by \cite[Theorem 3.7]{CCC} in terms of the characters of projective covers.  Consequently}, the character formulae of $T^{\mf b}_\la$ for $\mf{pe}(2)$ is minimal in the following sense: 	
 \begin{itemize}
 	\item[(1).] If $\mu$ is strongly linked to $\la$, then $(T^{\mf b}_\la:\nabla^{\mf b}_{\mu})>0$.  
	 \item[(2).]If $\langle\la,\alpha_1\rangle+1 =0$, then $(T^{\mf b}_\la:\nabla^{\mf b}_{\lambda-\ov \alpha_1})>0.$ 
	  \end{itemize} Here $(T^{\mf b}_\la: \nb^{\mf b}_\zeta)$ denotes the number of multiplicity of  $\nb^{\mf b}_\zeta$ in a parabolic (dual) Verma flag of $T^{\mf b}_\la$ for given $\la, \mu \in \h^\ast$, see Subsection \ref{subsect42}.

  As another application of Theorem A, we compute the  character formulae of  tilting modules in $\mc O^\fp$ for $\mf{pe}(3)$ (see Section \ref{Sect::tranProj}).  In particular,  
  the character formulae of tilting modules of non-$\mf b$-weakly-typical highest weights in $\mc O$ is minimal in the following sense:

 \begin{thmB}\label{thm::thmB} Let  $\la,\mu\in \h^\ast$, and $\beta, \gamma \in \Phi_\oa^+$. The following statements hold for $\mf g=\mf{pe}(3)$. 
 		
 1.      If $\mu\uparrow\la$, then $(T^{\mf b}_\la:\nabla^{\mf b}_{\mu})>0$.%
 	
 	2.     If $\langle\la, \beta\rangle+1 =0$, then $(T^{\mf b}_\la:\nabla^{\mf b}_{\lambda-\ov \beta})>0.$ 
 	
 	3.    If $\langle\la, \beta\rangle+1 =0$ and $w(\la -\ov \beta)\uparrow\la -\ov \beta$, then $( T^{\mf b}_{\lambda}:\nabla^{\mf b}_{w(\lambda-\ov \beta)})>0$.
 	
 		4.    If $\langle\la, \alpha_j\rangle+1 = \langle \la -\ov \alpha_j, \alpha_i \rangle+1 =0$ for some $1\leq  i<j\leq n-1$, then $(T^{\mf b}_\la:\nabla^{\mf b}_{\lambda -\ov \alpha_j -\ov \alpha_i})>0$.

 	5. If $\langle\la, \alpha_i \rangle+1 = \langle\la, \alpha_{i+1}\rangle+1 =0$ and $\beta:= \alpha_i+\alpha_{i+1}$, then we have $(T^{\mf b}_\la:\nabla^{\mf b}_{\lambda -\ov \alpha_{i} -\ov \beta})>0$ and $(T^{\mf b}_\la:\nabla^{\mf b}_{\lambda -\ov \alpha_{i+1} -\ov \beta})>0$.
 	
 	6. If $\langle \la, \alpha_i\rangle +1 =0$ for all $1\leq i\leq n-1$, then $(T^{\mf b}_{\la}:\nb^{\mf b}_{\la-(n-1)(\vare_1+\cdots +\vare_n)})=1.$ 
 \end{thmB}

 We refer the reader to \cite[Proposition 2.2]{ChWa18} for similar formulae for basic classical Lie superalgebras.  In fact, {\em 1} and {\em 6} in Theorem D are true for any $n$.  In addition, by Theorem A and Proposition \ref{prop::Jan}, one can deduce that {\em 2} and {\em 3} also hold for a specific class of simple roots $\beta$ for any $n$.
It is an interesting problem to generalize other statements in Theorem D, possibly with mild modifications, to any $n$.

 \subsection{}
 The paper is organized as follows. 
 
 In Section \ref{Sect::pre}, we set up notations and recall some preliminary results on the periplectic Lie superalgebras.  In Section \ref{Sect::parabolic}, we introduce the block decomposition in the full category $\mc O$.  Theorem \ref{thm::main} offers a   version for Theorem B.  Also, we reformulate Theorem A  as given in  Theorem \ref{thm24}.   Subsequently, the  proof of Theorem B is established in Subsection \ref{subSect::33}.

 In Section \ref{Sect::chwtw} some general consequences about the character of tilting modules of $\fp$-weakly-typical highest weights are given.  In particular, the proof of Theorem C is given in that of Theorem \ref{cor::ChP}.   Finally, the Sections \ref{Sect::tranProj} and \ref{Section::app}  are devoted to the proof of  Theorem D by providing the complete list of tilting characters in $\mc O^\fp$ for $\mf{pe}(3)$.

\vskip 0.5cm 
{\bf Acknowledgments.}  Both authors are grateful to Shun-Jen Cheng, Kevin Coulembier and Weiqiang Wang  for numerous interesting discussions. We have learned the validity of Proposition \ref{prop::Jan}  from Shun-Jen Cheng. We  thank Catharina Stroppel for useful remarks and comments on earlier versions of this paper. We also thank Li Luo, Bin Shu and East China Normal University for their hospitality during our visit in 2019 when a part of this project is carried out.  Chen is partially supported by MoST grant 108-2115-M-008-018-MY2.  Peng is partially supported by MoST grant  105-2628-M-008-004-MY4.

\section{Preliminaries}  \label{Sect::pre}
Throughout the paper the symbols  $\mathbb C$, $\mathbb Z$, $\mathbb N$, and  $\mathbb Z_{>0}:=\mathbb N\setminus\{0\}$ stand for the sets of all complex numbers, all integers,  non-negative and positive integers, respectively. Denote the abelian group of two elements by $\mathbb Z_2 =\{\oa,\ob\}$.  All vector spaces, algebras, tensor products, et
cetera, are over $\mathbb C$.

\subsection{Periplectic Lie superalgebra}\label{Intropen}
For positive integers $m,n$, the general linear Lie superalgebra $\mathfrak{gl}(m|n)$ may be realized as the space of $(m+n) \times (m+n)$ complex matrices
\begin{align*}
\left( \begin{array}{cc} A & B\\
C & D\\
\end{array} \right),
\end{align*}
where $A,B,C$ and $D$ are respectively $m\times m, m\times n, n\times m, n\times n$ matrices, with Lie bracket given by the super commutator.
Let $E_{ab}$ be the elementary matrix in $\mathfrak{gl}(m|n)$ with $(a,b)$-entry $1$ and other entries 0, for $1\leq a,b \leq m+n$.

The standard matrix realization of the periplectic Lie superalgebra $\pn$ inside $\mathfrak{gl}(n|n)$ is given
by
\begin{align}\label{plrealization}
 \mf g= \mf{pe}(n)=
\left\{ \left( \begin{array}{cc} A & B\\
C & -A^t\\
\end{array} \right)\| ~ A,B,C\in \C^{n\times n},~\text{$B^t=B$ and $C^t=-C$} \right\}.
\end{align}
Throughout the present paper, we fix the Cartan subalgebra $\mf h= \mf h_\oa \subset \mf g_\oo$ consisting of diagonal matrices. We denote the dual basis of $\mf h^*$ by $\{\vare_1, \vare_2, \ldots, \vare_n\}$ with respect to the standard basis
$$\{E_{ii}-E_{n+i,n+i}|~1\leq i \leq n \}\subset \pn.$$  The set of roots is given by
\begin{equation}\label{eqroots}
\Phi\;=\{\epsilon_i-\epsilon_j\,|\, 1\le i\not=j\le n\} \cup \{ \epsilon_i+\epsilon_j\,|\, 1\le i\le j\le n\} \cup  \{- \epsilon_i-\epsilon_j\,|\, 1\le i< j\le n\}. 
\end{equation} We recall the definitions \eqref{Roots1}, \eqref{Roots2} of the  subsets $\Phi_\oa$, $\Phi_\oa^+$ and $\Phi_\ob$. 
The Weyl group $W=\mf S_n$ of $\g$ is the symmetric group on $n$ symbols. Let $\rho:= \sum_{i=1}^n (n-i)\vare_i$ be the Weyl vector.
For any $\alpha\in\Phi_\oa$, we let $s_\alpha$ denote the corresponding reflection in $W$.  

We  fix a   Borel subalgebra $\fb_{\oa}$ of $\g_{\oa}=\mathfrak{gl}(n)$ consisting of matrices in \eqref{plrealization} with $B=C=0$ and $A$ upper triangular. Unless mentioned otherwise, all  parabolic subalgebras  are assumed to contain $\fb_{\oa}$.

Define the following subalgebras of $\g$:
\begin{align*}
\mf g_1:=
\{\begin{pmatrix}
0 & B \\
0 & 0
\end{pmatrix}|B^t=B\}\quad\mbox{and}\quad \mf g_{-1}:=
\{\begin{pmatrix}
0 & 0 \\
C & 0
\end{pmatrix}|C^t=-C\}.
\end{align*}
The standard Borel subalgebra $\mf b$ and the reverse Borel subalgebra $\mf b^r$ are defined by
$$\mf b:=\mf b_\oa \oplus \mf g_1,\qquad \mf b^r:=\mf b_\oa \oplus \mf g_{-1}.$$


A weight $\la\in \h^\ast$ is said to be dominant (resp. anti-dominant) if $\langle\la, \alpha \rangle \notin \mathbb Z_{<0}$ (resp. $\langle\la, \alpha \rangle \notin \mathbb Z_{>0}$), for all $\alpha \in \Phi_\oa^+$. Finally, we let $X=\sum_{i=1}^n\mZ\ep_i$,  and  $\omega_k:=\vare_1+\cdots +\vare_k$, for any $1\leq k\leq n$. 

\subsection{Representations}
\subsubsection{The BGG category}
Recall that we denote by $\cO=\cO(\fg,\fb)$ the BGG category of~$\fg$-modules which are finitely generated, semisimple as $\fh_{}$-modules and locally finite as $\fb_{}$-modules. This is thus the category of~$\fg$-modules which restrict to modules in the BGG category~$\cO^{\oa}:=\cO(\fg_{\oa},\fb_{\oa})$ of~\cite{BGG}.

  For any $\la \in \h^\ast$, we denote by $$M^{\oa}_\lambda=U(\fg_{\oa})\otimes_{U(\fb_{\oa})}\mC_{\la -\rho}\text{ and }\Delta_\la := U(\mf g)\otimes_{U(\mf b)}   \mC_{\la-\rho}\cong \text{Ind}_{\mf g_\oa\oplus \mf g_1}^{\mf g}M^\oa_\la,$$ the $\mf g_\oa$-Verma module and $\mf g$-Verma module of highest weight $\la -\rho$, respectively. Then the corresponding unique simple quotient of $M^\oa_\la$ and $\Delta_\la$ are respectively denoted by $L^\oa_\la$ and $L_\la$. 
Similarly,  we define the dual $\mf g$-Verma module by
$$\nabla_\la:= \text{Coind}_{\mf g_{\leq 0}}^{\mf g} ((M^\oa_\la)^{\vee}),$$  where $(M^\oa_\la)^{\vee}$ is the dual Verma module in~$\mc O^\oo$, see \cite[Subsection 3.2]{Hu08}.
Then we have $\Delta_\la \twoheadrightarrow L_\la$ and   $L_\la \hookrightarrow \nabla_\la$. 

The category $\mc O$ is a highest weight category (cf. \cite[Theorem 3.1]{CCC}) with standard objects $\Delta_\la$ and costandard objects $\nb_\la$. However,  $\mc O$ does not admit a simple-preserving duality by \cite[Corollary 4.6]{CCC}. It was shown in \cite[Proposition 3.4]{CCC} that the costandard objects $\nb_\la$ are still images of standard objects for different highest weight structure of $\mc O$ under a natural duality.

For a given object $V \in \mc O^\oa$, we extend $V$ trivially to a $(\mf g_\oa \oplus \mf g_1)$-module  and define the {\em Kac module} $$K(V):=\Ind_{\mf g_\oa\oplus \mf g_1}^{\mf g}V.$$ This defines an exact functor $K(\cdot): \mc O^\oa \rightarrow \mc O$ which is called the Kac functor.  In particular, we write  $K_\la := K(L^\oa_{\la}),$ for any $\la \in \h^\ast$. 
 
Denote the irreducible module of $\mf b^r$-highest weight $\la -\rho$ by $L_\la^{\mf b^r}$. It is well-known that the conditions for $\la$ and $\mu$ satisfying  $L_\la^{\mf b^r} = L_\mu$ can be explicitly described by an algorithm involving odd reflections. Throughout the paper, we refer to \cite[Section 2.2]{PS89} for a treatment of odd reflections for  Lie superalgebras. Particularly, in \cite[Lemma 1]{PS89} the effect on the highest weight of a simple $\pn$-module under
odd reflection and inclusion was computed.

 Also, for any $\la \in \h^\ast$  by \cite[Lemma 5.2]{CC} we have   
 \begin{align}
 &\ch \nabla_\la =\sum_{\kappa \in S} \ch \Delta_{\la -\kappa}, \label{Eq::Delnb}
 \end{align}  here $S:=\{\kappa = \sum_{i=1}^{n} \kappa_i \vare_i~|\kappa_i \in \{0,2\} \text{ for all }i\}$.

  Finally, for a given module $M \in \mc O$ we denote its socle and radical by ${\rm soc}M$ and ${\rm rad} M$, respectively.  We denote by $\mc O_{\mathbb Z}$ the full subcategory of $\mc O$ consisting of modules with weights in $X$.

\subsubsection{$\mathbb Z$-gradation} We introduce a natural $\mathbb Z$-grading for objects in the category $\mc O$. For each $\mu = \sum_{i=1}^n \mu_i \vare_i\in \h^\ast$ we set $|\mu|$ to be the sum of all components of the weight $\mu -\rho$, that is, $$|\mu|:=\mu_1+\mu_2+\cdots +\mu_n -\sum_{i=1}^n (n-i).$$ For any $M\in \mc O$, we let $M=\bigoplus_{z\in \C} M_{[z]}$ denote the eigenspace decomposition of $M$ with respect to the grading operator $d:=   \sum_{i=1}^n E_{i,i}-E_{n+i,n+i}\in \mf g_\oa$. Namely, $$M_{[z]}:=\{m\in M|~ dm =zm\}.$$ Note that  $M_{[z]}\in \mc O^\oo$ for any $z\in \C$.  
In particular,  we have $$K_\mu=\bigoplus_{k=0}^{\frac{n(n-1)}{2}}(K_{\mu})_{[|\mu| -2k]},~\Delta_\mu=\bigoplus_{k=0}^{\frac{n(n-1)}{2}}(\Delta_{\mu})_{[|\mu| -2k]},~\nabla_\mu=\bigoplus_{k=0}^{\frac{n(n+1)}{2}}(\nabla_{\mu})_{[|\mu| -2k]}$$ 

Note that $L_\la$ is $\mathbb Z$-graded and the number of its $d$-eigenvalues is equal or less than $\frac{n(n-1)}{2}+1,$ for any $\la \in \h^\ast.$  
Moreover, $L_\la=(L_\la)_{[|\la|]}\oplus (L_\la)_{[|\la|-2]}\oplus  \cdots \oplus (L_\la)_{[m]},$ for some $m$. Here 
$(L_\la)_{[|\la|]} \cong L^\oa_\la$ and   $(L_\la)_{[m]}$ are irreducible $\mf g_\oa$-modules.



\section{Blocks in   $\mc O^{\mf p}$} \label{Sect::parabolic} In this section, we investigate the blocks in the  parabolic version of category $\mc O$ for $\pn$. 
From now on, we set $\mf p$ to be a {\em reduced} parabolic subalgebra as defined in \cite[Subsection 1.4]{CCC} with the purely even Levi subalgebra $\mf l = \mf l_\oa\cong \bigoplus_{i=1}^{k}\gl(n_i)$ for some positive integers $k$ and $n_i$'s. Namely, if we set $n_0$ to be zero then we have $$\Phi^+(\mf l)=\{\vare_i -\vare_j|~n_0+n_1+\cdots + n_q+1 \leq i<j \leq n_0+n_1+\cdots+ n_{q+1},\text{ for some $0\leq q\leq k-1$}\}.$$ The {\em parabolic category} $\mc O^\fp  = \mc O(\mf g, \mf p)$ is the full subcategory of $\mc O$ consisting of finitely generated $\mf g$-modules on which $\fp$ acts locally finitely. Without loss of generality, we assume that $\mf p_\ob=\mf  g_1$. 
 
\subsection{Parabolic category $\mc O^\fp$} \label{Sbusect::parabolicCate}

For any $\la \in \h^\ast$, let $L^{\mf l}_\la$ denote the irreducible $\mf l$-module with highest weight $\la -\rho$. 
	Recall the set of {\em $\fp$-dominant weights} defined in \cite[Section 3]{CCC}:
	$$\Sigma^+_\fp:=\{\la\in \mf h^\ast|~ \langle\la, \alpha\rangle \in \mathbb Z_{>0},\text{ for all } \alpha \in \Phi^+(\mf l) \}.$$ Observe that $\la \in \Sigma_\fp^+  \Leftrightarrow \text{dim} L_\la^{\mf l} <\infty.$ The corresponding  parabolic category $\mc O^\fp$ is the Serre subcategory of $\mc O$ generated by $\{L_\la \, | \,  \la \in \Sigma_\fp^+\}$. Also, we define $\mc O^\fp_{\Z}$ to be the Serre subcategory of $\mc O$ generated by  $\{ L_\la \, |\,  \la \in \Sigma_\fp^+ \cap X\}$.

For $\la \in \Sigma_\fp^+$, we define respectively the parabolic Verma $\mf g_\oa$-module and  the parabolic Verma $\mf g$-module by 
\begin{align}
&M^{\fp_\oa}_\la : =\Ind_{\mf p_\oa}^{\mf g_\oa}L^{\mf l}_\la, ~\Delta^\fp_\la:=\Ind_{\mf p}^{\mf g}L^{\mf l}_\la. \label{Eq::paraV}
\end{align} We may note that  $\Delta^\fp_\la  \cong \Ind_{\mf g_\oa+\mf g _1}^{\mf g}M^{\fp_\oa}_\la$. Also, we define the dual parabolic Verma module by $$\nabla^\fp_\la:=\Ind_{\mf g_\oa+\mf g _{-1}}^{\mf g}((M^{\fp_\oa}_\la)^{\vee}\otimes \Lambda^{\text{top}}\mf g_1^\ast)\cong \Coind_{\mf p_\oa+\mf g_{-1}}^{\mf g}L^{\mf l}_\la.$$ We have $\Delta^\fp_\la \twoheadrightarrow L_\la$ and   $L_\la \hookrightarrow \nabla^\fp_\la$. 
If $\la\in \Sigma^+_\fp$, we let $\mc O^\fp_\la$ denote the (indecomposable) block in $\mc O^\fp$ containing $L_\la$. 
We denote by $\mc O_\la:=\mc O_\la^{\mf  b}$ the block  in the full category $\mc O$ containing $L_\la$.

\subsection{A description of blocks in $\mc O^\fp$}   \label{subsect::link}   We recall the equivalence relation $\sim$ on $\h^*$ defined in~\cite[Subsection 5.2]{CC} which is transitively generated by 
$$\begin{cases}\la\sim \la \pm2\vare_k, &\mbox{ for~$1\le k\le n$;}\\
\la\sim w  \la,&\mbox{ for~$w\in W^{[\lambda]}$.}
\end{cases}$$
Here $W^{[\lambda]}$ denotes the integral Weyl group associated to $\la$. 

 We set $\partial_0:=0$, and for $1\leq i\leq n$ we let 
$$\partial_i:= \sum_{j=1}^{i} \vare_j$$ 
be  weights of $\pn$. Let $\mc C_n$ denote the set  $\{\partial_{0},~ \partial_1,~\partial_{2}, ~\ldots, ~ \partial_{n}\}$ for $\pn$. For a given weight $\la\in X$, we let $N_{\rm odd}(\la)$ be the number of odd integers in $\{\la_1, \la_2,\ldots, \la_n\}$. In particular, we have $N_{\rm odd}(\partial_i) =i,$ for $0\leq i\leq n.$

 The following lemma establishes the block decomposition of the full category $\mc O$ and shows that the blocks of $\mc O_\Z$  are indexed by the  $\partial_0,\partial_1,\ldots, \partial_n$.  \begin{lem}{\em (}\cite[Theorem 5.4]{CC}{\em )} \label{lem::CC54} Let $\la, \mu \in \h^\ast$. 
		Then $L_\la$ and $L_\mu$ are in the same block of $\mc O$ if and only if $\lambda\sim\mu$. 
		
		In particular,  the category $\mc O_\Z$ has exactly $n+1$ blocks:
		\begin{align}
		&\mc O_\Z =\bigoplus_{i=0}^n \mc O_{\partial_i}.
		\end{align}
		Namely, we have  $L_\la\in \mc O_\mu \Leftrightarrow N_{\rm odd}(\la)=N_{odd}(\mu),$ for any  $\la,\mu\in X$.
	\end{lem}  
\begin{proof} The linkage $L_\la \sim L_\mu \Leftrightarrow \la\sim \mu$ is given by \cite[Theorem 5.4]{CC}. It remains to show the claims for block decomposition of $\mc O_\Z$. 
	
	Let $\partial'_0=\rho$. For $1\leq i\leq n$ and $\partial'_i:= \big( \sum_{j=1}^{i} (i+1-j)\vare_j  \big)+\rho$
	be (shifted) 2-cores for $\pn$ appearing in \cite[Section 7.1]{Co}. By \cite[Theorem 5.4]{CC} we have block decomposition:
	\begin{align}
	&\mc O_\Z =\bigoplus_{i=0}^n \mc O_{\partial'_i}.
	\end{align}
	Using the linkages $\la \sim \la\pm 2\vare_q$ for $1\leq q\leq n$, one can deduce a  permutation $\sigma$ on the set $\{0,1,\ldots, n\}$ such that $\mc O_{\partial'_i} =\mc O_{\partial_{\sigma(i)}},$ for $1\leq i\leq n.$  This completes the proof. 
	\end{proof}
 For $\la \in \Sigma_\fp^+$, we define $$\g^{[\la]}_\oa = \mf h\oplus \bigoplus\limits_{\substack{\alpha\in \Phi_\oa,  \\  \langle\la, \alpha\rangle\in \Z}}\mf g_\alpha.$$ 
	Let $\mc O_{\la+X}$ be the Serre subcategory of $\mc O$ generated by $L_{\mu}$ for $\mu \in \la+X$. Then Lemma \ref{lem::CC54}  gives rise to  the following block decomposition \begin{align} 
	&\mc O_{\la+X} =\bigoplus_{\mu \in \la+X/ \sim} \mc O_{\mu}. \label{Eq::blocks1}
	\end{align}

To further describe the block decomposition (\ref{Eq::blocks1}), we  assume that $\mf g^{[\la]}_\oa$ is a Levi subalgebra of $\g_\oa$ with 
$$\mf g^{[\la]}_\oa\cong \bigoplus_{i=1}^{k}\gl(n_i).$$ 
Let $r_1,\ldots, r_k \in \C$ be given as follows: 
$$r_1=\la_1,~r_2 = \la_{n_1+1},~ r_{3}= \la_{n_1+n_2+1},~ \ldots,~ r_k=\la_{n_1+n_2+\cdots +n_{k-1}},$$ which are all  representatives in the conjugacy classes  $\{\la_i|~1\leq i\leq n\}/ \Z$. For given numbers $r\in \C$, $b\in \Z_>0{}$ and  a given weight $\mu=\mu_1\vare_1+\mu_2\vare_2+\cdots + \mu_a\vare_a$ with $a+b \leq n$, we define $$\mu[r;b]:=(r+\mu_1)\vare_{1+b} +(r+\mu_2)\vare_{2+b}+\cdots +(r+\mu_{a})\vare_{a+b}.$$ With these fixed numbers $r_1,\ldots, r_k$, we now  identify elements in $\mc C_{n_1}\times\mc C_{n_2}\times \cdots \times \mc C_{n_k}$ with the following specific weights in $\la+X$: $$(\partial_{i_1}, \partial_{i_2},\ldots, \partial_{i_k}):= \partial_{i_1}[r_1;0]+\partial_{i_2}[r_2;n_1]+\partial_{i_3}[r_3;n_1+n_2]+\cdots +\partial_{i_k}[r_k;n_1+n_2+\cdots +n_{k-1}],$$
	where $\partial_{i_j}\in \mc C_{n_j}$ for $1\leq j\leq k.$ By Lemma \ref{lem::CC54}, a description of  \eqref{Eq::blocks1} is then given  as follows, \begin{align}
	&\mc O_{\la+X} =\bigoplus_{\partial\in  \mc C_{n_1}\times \cdots \times \mc C_{n_k}} \mc O_{\partial},  \label{Eq::blocks2}
	\end{align} where $\mc O_{\partial}$ denotes the block  in $\mc O$ containing $L_\partial$. Similarly, we define $\mc O_{\partial}^\fp$ to be the block  in $\mc O^\fp$ containing $L_\partial$.

In this section, we shall establish the following theorem. Its proof will be given in the next subsection.
\begin{thm}\label{thm::main}
	Let $\la \in \Sigma_\fp^+$. Suppose that $\mf g^{[\la]}_\oa$ is a Levi subalgebra of $\g_\oa$ with $\mf g^{[\la]}_\oa\cong \bigoplus_{i=1}^{k}\gl(n_i)$. Then the mapping 
		$\partial \mapsto \mc O^\fp_\partial$ gives rise to a bijection from $\mc C_{n_1}\times  \cdots \times \mc C_{n_k}$ to the set of blocks in $\mc O_{\la+X}^\fp$. 
		Namely, each block in $\mc O_{\la+X}$ remains indecomposable when restricted to   $\mc O^\fp_{\la+X}$. 
                 In particular, the number of blocks in  $\mc O^\fp_{\la+X}$ is given by
                 $$\prod_{i=1}^{k} (n_i+1).$$
\end{thm} 
\vskip 0.3cm

\subsection{Proof of Theorem \ref{thm::main} and Theorem B} \label{subSect::33} This subsection is devoted to the proof of Theorem \ref{thm::main} and Theorem B. The following lemma is useful.

\begin{lem} \label{lem21}	Let $M\in \mc O$. For any $\la, \mu \in \h^\ast$ we have the following facts.
	\begin{enumerate} 
		\item[(1).]  If $M_{[z+1]} =0$  and $L^\oa_\la$ is a summand of $\emph{soc}_{\mf g_\oa}M_{[z]}$, then $\emph{Hom}_{\mc O}(K_\la, M)\neq 0$.  
		\item[(2).]   If $M =M_{[z]}\bigoplus \left( \bigoplus_{z>c} M_{[c]}\right)$ with $[M_{[z]}:L^\oa_\mu]\neq 0$, then $[M:L_\mu]\neq 0.$
	\end{enumerate}
\end{lem}
\begin{proof} 
The first statement follows easily from adjunction.
To prove the second statement, we use induction on the length $\ell(M_{[z]})$ of  the $\mf g_\oa$-composition series of $M_{[z]}$. 
The initial case $\ell(M_{[z]}) =1$ follows from {\em (1)}.	
	
	Assume that  $L^\oa_\la$ is a summand of $\text{soc}_{\mf g_\oa}M_{[z]}$ with $\la \neq \mu$. By \emph{(1)}, we have an epimorphism $K_\la \twoheadrightarrow N$ and a short exact sequence  in $\mc O$:
	$$0\rightarrow N \rightarrow M \rightarrow Q\rightarrow 0.$$ 
	It leads to the existence of the short exact sequence in $\mc O^\oa$:
	$$0\rightarrow L^\oa_\la\rightarrow M_{[z]} \rightarrow Q_{[z]} \rightarrow 0.$$
	Now $\ell(Q_{[z]}) = \ell(M_{[z]})-1$. By induction, if $\mu\in \h^\ast$ satisfies  $[Q_{[z]}:L^{\oa}_\mu]\neq 0$ then  $[Q:L_\mu]\neq 0$.  This completes the proof. 
\end{proof}

For $X, Y \in \mc O^\oa$, we define the following ordering: 
\begin{align} \label{Eq::chord}
&\ch X > \ch Y \text{ if } \ch X =\ch Y +\ch Z, \text{ for some } Z\in \mc O^\oa.
\end{align} 
Similarly, we define $\ch X - \ch Y > \ch X' -\ch Y'$ if $\ch X -\ch Y =  \ch X' -\ch Y' +\ch Z$, for some $Z\in \mc O^\oa.$

The following theorem and its proof imply  the validity of Theorem A in Section \ref{Sect::Introduction}.

\begin{thm}  \label{thm24} Let  $\la \in \Sigma^+_\fp$. For any $1\leq q\leq n$, define $A_q:=\{q\leq j\leq n|~\la_q =\la_j \}$. Then we have the following facts.
	\begin{itemize}
		\item[(i)] Let $1\leq q\leq n$ be given such that   $\la -2\vare_q\in \Sigma_\fp^+$. If $\langle\la,\alpha_j\rangle\neq 1$ for any $j\in A_q$, then    
		$$[\nb^\fp_{\la}:L_{\la-2\vare_q}]>0.$$ 
		In particular, we have $[\nb_{\la}:L_{\la-2\vare_q}]>0.$
			\item[(ii)]  Let $\alpha_i\in \Phi^+(\mf l)$ be given such that $\la -\ov \alpha_i\in \Sigma^+_\fp$ and  $\langle \la, \alpha_i\rangle =1$. If $A_i =\{i\}$, then 
			$$ [\Delta^\fp_{\la}:L_{\la-\ov \alpha_i}]>0.$$  
			In particular, we have $[\Delta_{\la}:L_{\la-\ov \alpha_i}]>0.$
	\end{itemize}
\end{thm}
\begin{proof}
	We prove (i) in detail here where (ii) can be proved in a similar way. Note that $$(\Delta^\fp_\la)_{[|\la|-2]} = \Lambda^1 \mf g_{-1}\otimes M^{\fp_\oa}_\la,~(\nabla^\fp_\la)_{[|\la|-2]} \cong \Lambda^1 \mf g_1^\ast\otimes M^{\fp_\oa}_\la.$$ We shall proceed with some estimates of  characters   of summands in these layers.
	

	Set $A^\la =\bigoplus_{\mu\in \h^\ast}[M^{\fp_\oa}_\la:L^\oa_\mu] L_\mu$. By Lemma \ref{lem21}, we have $A^\la \in \mc O^\fp$ such that $${\ch}\nabla^\fp_\la>\ch A^\la,~\ch\Delta^\fp_\la > \ch A^\la, $$   $[A^\la:L_\la]\neq 0$ and   $\ch(A^\la)_{[|\la|]}= \ch M^{\fp_\oa}_\la$.
	 It follows that 
	\begin{align*}
	&\text{ch}({\nabla^\fp_\la})_{[|\la|-2]}-\text{ch}({A^\la})_{[|\la|-2]}> \text{ch}\Lambda^1 \mf g_1^\ast\otimes M^{\fp_\oa}_\la - \text{ch}\Lambda^1 \mf g_{-1}\otimes M^{\fp_\oa}_\la.
	\end{align*}
	We note that as $\mf l$-decompositions there is an $\mf l$-module $C$ such that  $$\Lambda^1 \mf g_1^\ast = \bigoplus_{i=1}^{k}S^2(\C^{n_i\ast}) \oplus C, ~\Lambda^1 \mf g_{-1}=\bigoplus_{i=1}^{k}\Lambda^2(\C^{n_i\ast})\oplus C.$$
	By Pieri's rule (cf. \cite[Formula (5.16)]{Mac}), it follows that 
	\begin{align*}
	& \text{ch}\Lambda^1 \mf g_1^\ast\otimes L^{\mf l}_\la - \text{ch}\Lambda^1 \mf g_{-1}\otimes L^{\mf l}_\la\\
	&=\text{ch}(\bigoplus_{i=1}^{k}S^2(\C^{n_i\ast})\otimes L^{\mf l}_\la) -\text{ch}(\bigoplus_{i=1}^{k}\Lambda^2(\C^{n_i\ast})\otimes L^{\mf l}_\la)\\
	&=\sum_{\la-2\vare_i\in \Sigma^+_\fp} \text{ch}L^{\mf l}_{\la-2\vare_i} - \sum_{ \substack{\alpha_i \in \Phi^+(\mf l), \\\langle\la, \alpha_i \rangle=1,~\la-\ov \alpha_i\in \Sigma^+_\fp}}\text{ch}L^{\mf l}_{\la-\ov \alpha_i}.
	\end{align*}
	
	As a consequence, we have 
	\begin{align} \label{Eq::2020}
	&\text{ch}({\nabla^\fp_\la})_{[|\la|-2]}-\text{ch}({A^\la})_{[|\la|-2]} > \sum_{\la-2\vare_i\in \Sigma^+_\fp}  \text{ch} M^{\fp_\oa}_{\la-2\vare_i} -  \sum_{\substack{ \alpha_i \in \Phi^+(\mf l), \\ \langle\la, \alpha_i \rangle=1,~\la-\ov \alpha_i\in \Sigma^+_\fp}}\text{ch}M^{\fp_\oa}_{\la-\ov \alpha_i}. 
	\end{align}

	With some straightforward computation, it is not hard to show that the assumption of (i) is equivalent to the following statement:
	\begin{center}
	$\la-2\vare_q$ is not  strongly liked to  $\la-\ov \alpha_i$ for any $1\leq i\leq n.$ 
	\end{center}
	By the linkage principle for $\mc O^\oa$, we have  $[M^{\fp_\oa}_{\la-\ov \alpha_i}: L^\oa_{\la-2\vare_q}]=0$ for any $1\leq i\leq n-1$. This implies that the coefficient of $\text{ch}L^{\oa}_{\la-2\vare_q}$ in the  expression of  $\text{ch}(\nabla^\fp_\la)_{[|\la|-2]} - \text{ch}(A^\la)_{[|\la|-2]}$ is positive. 
	
	By our construction $\text{ch}\nabla^\fp_\la -\text{ch}A^\la$ is a character of a $\mf g$-module. Since  $(\nabla^\fp_\la)_{[|\la|]} = (A^\la)_{[|\la|]}$ and $(\nabla^\fp_\la)_{[|\la|-1]} = (A^\la)_{[|\la|-1]}=0$, we have $M = M_{[|\la|-2]}\bigoplus \left( \bigoplus_{|\la|-2>c}M_{[c]}\right) $. Together with the fact that $$[M_{[|\la|-2]}: L^{\oa}_{\la-2\vare_q}]> 0,$$  we may conclude from Lemma \ref{lem21}  that $[M:  L_{\la-2\vare_q}]> 0.$ Consequently, we have   $[\nabla^\fp_\la:  L_{\la-2\vare_q}]> 0$. This completes the proof. 
			
\end{proof}


We now in a position to prove Theorem \ref{thm::main}. Recall that we assume $\la \in \Sigma_\fp^+$ and $\mf g^{[\la]}_\oa$ is a Levi subalgebra of $\g_\oa$ with $\mf g^{[\la]}_\oa\cong \bigoplus_{i=1}^{k}\gl(n_i)$. Also note that $\mf l \subseteq \mf g_\oa^{[\la]}$.

\begin{proof}[Proof of Theorem \ref{thm::main}]
	Let $\partial,~\partial' \in \mc C_{n_1} \times \cdots \times \mc C_{n_k}$. If  $\partial\neq \partial'$ then $\mc O_{\partial}^\fp\neq \mc O_{\partial'}^\fp$ by    the decomposition \eqref{Eq::blocks2}. 
		In particular, the number of blocks in $\mc O^\fp_{\la+X}$ is at least $\prod_{i=1}^{k} (n_i+1)$. It remains to show that each block in $\mc O_{\la+X}^\fp$ is of the form $\mc O_\partial^\fp$, for some $\partial\in  \mc C_{n_1} \times \cdots \times \mc C_{n_k}$.

	We first assume that $\mf g^{[\la]}_\oa =\mf g_\oa$, that is, $\la\in X$. We claim that in this case  the block $\mc O_\la^\fp$ contains  $L_{\partial_i}$, for some $0\leq i\leq n$.  	To see this, we first observe that by Theorem \ref{thm24} the simple module  $L_{\la-2N\vare_n}$ lies in  $\mc O_\la^\fp$ for any positive integer $N$. Therefore, there exists some $N_1\in \Z$ such that $L_{\la-2N_1\vare_n} \in \mc O_\la^\fp$ with $\langle \la-2N_1\vare_n, \vare_i-\vare_{n}\rangle>0$ for any $i<n$.  For a fixed $1\leq q\leq n$,  if the weight $\la'\in \h^\ast$ satisfies the condition $\la'_q>\la'_j$ for all $j>q$, then $L_{\la'-2\vare_q}\in \mc O^\fp_{\la'}$ by Theorem \ref{thm24} (ii).
		Repeatedly applying the argument above to $q=n-1,n-2,\ldots,1$, we eventually find  a simple module $L_{\mu} \in \mc O_\la^\fp$ such that $\langle\mu,\alpha_i\rangle \in \Z_{> 0},$ for all $1\leq i\leq n-1$.  Consequently, we may conclude from Lemma \ref{lem::CC54} that  $L_{\partial_i} \in \mc O_\la^\fp$ for some $0\leq i\leq n$, as desired.
	
		We now assume that $\la \notin X$. Our goal is to show that $\mc O_\la^\fp = \mc O_\partial^\fp$ for some $\partial\in \mc C_{n_1}\times \cdots \times\mc C_{n_k}.$ Note that $W^{[\la]}\cong \mf S_{n_1} \times\mf S_{n_2}\times\cdots \times\mf S_{n_k}$.	Using a similar argument as above, we obtain a  simple module $L_{\mu} \in \mc O_\la^\fp$ such that $\langle\mu,\alpha_i\rangle \not\in \Z_{\leq 0},$ for all $1\leq i\leq n-1$. Therefore we may conclude from Lemma \ref{lem::CC54}  that there exists $\partial\in \mc C_{n_1}\times \cdots \times \mc C_{n_k}$ such that $L_{\partial} \in \mc O_\la^\fp$. This completes the proof.
\end{proof}

	\begin{proof}[Proof of Theorem B] Let $\mu \in \h^\ast$ and $1\leq i\leq n-1$ such that $\langle\mu,\alpha_i\rangle\notin \Z$. 
	By an argument very similar to that employed in the proof of \cite[Proposition 4.4]{CC} (also, see  \cite[Section 3.6]{CMW}), one can show that there is an equivalence $\mc O_\mu \cong \mc O_{s_{\alpha_i}\mu}$ sending $L_\zeta$ to $L_{s_{\alpha_i}\zeta}$, for any $L_\zeta \in \mc O_\mu$. As a consequence, there exists   $\la \in \h^\ast$ such that $\mf g^{[\la]}_\oa$ is a Levi subalgebra of $\g_\oa$ and  $\mc O_\mu\cong \mc O_\la$ which restricts to  equivalences of parabolic subcategories of $\mc O_\mu$ and $\mc O_\la$. Consequently, the proof   follows from Theorem~\ref{thm::main}.
\end{proof}

	\begin{rem}The following facts generalizing  Brundan's indecomposablity theorem  for $\mf g_\oa$  provide   an alternative approach to  Theorem B. \\
		(1). By \cite{Br4} (see also \cite[Section 5.1]{CoM2}), every integral block of $\mc O^\oa$ remains indecomposable  when restricted to parabolic subcategory of $\mf g_\oa$-modules (whenever it is non-zero). Using reduction procedure developed  in \cite[Proposition 2.3]{CMW}, one   generalize  to blocks of modules of non-integral weights. \\
		(2). By \cite[Theorem 49]{CCM}, we may conclude that the Kac functor give rise to an  embedding $$K(-): {\rm Ext}^1_{\mc O^\oa}(L^\oa_\la,L^\oa_\mu) \hookrightarrow {\rm Ext}^1_{\mc O^\fp}(L_\la,L_\mu),$$ for any $\la,\mu \in \Sigma_\fp^+.$ 
			An alternative proof of Theorem B then follows from Theorem \ref{thm24}-(i). 
\end{rem}


\subsection{A linkage principal in $\mc O$}

We will describe a linkage principal about composition factors in Verma modules, which will be  helpful in computing characters in sequel. 
\begin{cor} \label{cor::oddlinkage}   Let $\mf g =\mf{pe}(3)$. If $\la \in \h^\ast$ satisfies $\langle\la, \alpha_i\rangle =1$, then $[\Delta_{\la}:L_{\la-\ov \alpha_i}]>0$. 
\end{cor}
\begin{proof}  Let $\la =\la_1\vare_1+\la_2\vare_2 +\la_3\vare_3$.	If either $i=2$ or $i=1$ with $\la_3\neq \la_1$, then we have 
	$[\Delta_{\la}: L_{\la - \ov \alpha_i}]>0$ by  Theorem \ref{thm24}. 
	
	 We now assume that $i=1$ with $\la_1=\la_3$. By using the algorithm of odd reflection given in \cite[Section 2.2]{PS89}, we found that 
	$$L_{\la_1\vare_1+(\la_1-1)\vare_2+(\la_1+1)\vare_3} = \text{soc}K_{(\la_1+1)\vare_1+\la_1\vare_2+(\la_1+1)\vare_3}.$$  This completes the proof.   
\end{proof}


\begin{prop} \label{prop::Jan}    Let $\la, \mu \in \h^\ast$  
be such that $\mu  \uparrow  \la$. If  $[\Delta_\mu: L_\eta] >0$ for some  $\eta\in \h^\ast$, then $[\Delta_\la :L_\eta] >0.$ 
\end{prop}
\begin{proof}
	Consider the Jantzen filtration of $M^\oa_\la$ as given in \cite[Section 5.3]{Hu08}
		\begin{align}\label{Eq::g0Jan}
	& M^\oa_\la =: M^0\supset M^1 \supset M^2 \supset \cdots,
	\end{align} with $M^1 =\text{rad}M^\oa_\la$ and 
	\begin{align}
	&\sum_{i>0} \ch M^i =\sum_{\la <s_\alpha\la} \ch M^\oa_{s_\alpha \la}, 
	\end{align}where $\alpha$ takes over positive even roots.
	Applying the Kac functor to (\ref{Eq::g0Jan}), we have the following filtration
	\begin{align}
	&\Delta_\la \supset K(M^1) \supset K(M^2) \supset \cdots.
	\end{align} We obtain the following sum formula 
	\begin{align}
	&\sum_{i>0} \ch K(M^i) =\sum_{\la <s_\alpha\la} \ch \Delta_{s_\alpha \la}, 
	\end{align}where $\alpha$ takes over positive even roots. This completes the proof. 
\end{proof}

\section{Characters of $\mf p$-weakly-typical  tilting modules} \label{Sect::chwtw}

\subsection{$\mf p$-weakly-typical weights}\label{subsubSect::wtw}
\subsubsection{A characterization}

It is known in \cite[Lemma 5.11]{CC} (see also \cite[Lemma 3.1]{Se02}) that 
if $\la$ is $\mf g_\oa$-weakly-typical, then $K_{\la} = L_{\la}$. We now prove the converse statement.  

\begin{prop} \label{lem::wt}  Let $\la \in \Sigma^+_\fp$. Then the following are equivalent:
	\begin{itemize}
		\item[(1).] $\la$ is $\mf p$-weakly-typical.
		\item[(2).] $-w_0^\fp\la$ is    $\mf g_\oa$-weakly-typical.
		\item[(3).] $K_{-w_0^\fp\la} = L_{-w_0^\fp\la}$.
	\end{itemize}
\end{prop}
\begin{proof} (1) and (2) are equivalent by definition. It remains to show that (3) implies (2).
	
		Assume that  $K_{-w_0^\fp\la}$ is simple. Then the $\mf b^r$-highest weight of $L_{-w_0^\fp\la}$ is $-w_0^\fp\la-(n-1)\omega_n$.
	By straightforward computation using the odd reflection given in  \cite[Section 2.2]{PS89}, one deduces that  $-w_0^\fp\la$ is $\mf g_\oa$-weakly-typical. This completes the proof. 
\end{proof}

\subsubsection{Simplicity of parabolic Verma module} Let $\la\in \h^\ast$ be an integral weight. Using the Kac functor, we observe that $\text{soc}\Delta_\la =L_\mu$ where  $\mu \in W\la$ is anti-dominant. This implies that every monomorphism between Verma modules is an isomorphism, and so $\Delta_\la$ is irreducible if and only if $\la$ is  anti-dominant. 
The following corollary  reduces the simplicity of parabolic Verma modules of $\fp(n)$ to that of $\gl(n)$, and the later  can be completely  described by Jantzen's simplicity criterion for parabolic Verma modules over $\mf g_\oa$ in \cite[Satz 4]{Ja}.
\begin{cor}
	Let $\la \in \Sigma^+_\fp$. Then $\Delta^\fp_\la$ is irreducible if and only if $\la$ is $\mf g_\oa$-weakly-typical and $M_\la^{\fp_\oa}$ is irreducible.  
\end{cor}
\begin{proof}
	Applying the Kac functor, we have \[[\Delta_\la^\fp:L_\mu] =\sum_{\zeta\in \h^\ast}[M_{\la}^{\fp_\oa}: L^\oa_\zeta][K_\zeta:L_\mu],\] for any $\mu \in \h^\ast.$ The proof follows from Proposition \ref{lem::wt}.
\end{proof}

\subsubsection{Irreducible modules of $\mf g_\oa$-weakly-typical highest weights}
 \begin{lem} \label{prop::tyKac}
	Let $\la\in \h^\ast$ be $\mf g_\oa$-weakly-typical. Then $[K_\mu:L_{\la}]\neq 0$ if and only if $\mu =\la$. In particular, we have $[\Delta_\mu : L_{\la}] = [M^\oa_\mu: L^\oa_{\la}]$, for any $\mu \in \h^\ast$.  
\end{lem}
\begin{proof} Assume that  $[K_{\mu}:L_{\la}]>0$. 
	Then $L_{\la}$ is a subquotient of $K_\mu$ and we have a short exact sequence in $\mc O$ $$0\rightarrow N \rightarrow M \rightarrow L_{\la} \rightarrow 0,$$
	where $M$ is a submodule of $K_\mu$. Here we may identify $L_{\la}$ as a  semisimple $d$-submodule of $K_\mu$. 
	
	By Proposition \ref{lem::wt}, we have $K_{\la}=L_{\la}$, which implies that the number of eigenvalues of $d$ acting on $L_{\la}$ is $\frac{n(n-1)}{2}+1$. These eigenvalues are explicitly given by
	$$|\mu|, ~|\mu|-2, \ldots,~ |\mu|-n(n-1).$$  
	Consequently, we have  $$(L_{\la})_{[|\mu|-n(n-1)]}\cap \text{soc} K_{\mu} \neq 0 \qquad \text{and} \qquad (L_{\la})_{[|\mu|]}\cap  (K_\mu)_{[|\mu|]} \neq 0.$$ 
	The former implies that  $N =0$ and hence $M =L_{\la}$  since $\text{soc}K_\mu = U(\mf g)\cdot (K_\mu)_{[|\mu|-n(n-1)]}$ is irreducible. 
	 The latter implies     $K_\mu = L_{\la}$, as desired. Therefore we have $\mu=\la$. This completes the proof. 
\end{proof}

The following corollary is a direct consequence of Lemma \ref{prop::tyKac}.
\begin{cor}
	$\Delta_{\la} = \sum_{\mu \uparrow \la} [M^\oa_\la: L^\oa_{\mu}] L_{\mu}$  for any typical weight $\la\in \h^\ast$.  
\end{cor}
 
\vskip 0.5cm


\subsection{Tilting modules of $\mf p$-weakly-typical highest weight}  \label{subsect42}

Let $M\in \mc O^\fp$.  
We say $M$ has a {\em parabolic (dual) Verma flag} 
if it has a filtration by submodules
$$0=M_0\subset M_1\subset\cdots\subset M_k=M,$$
where $M_i/M_{i-1}$ is a parabolic (dual) Verma module for each $1\le i\le k$. We denote by~$(M:\Delta^\fp_\lambda)$ and $(M:\nb^\fp_\lambda)$ the number of indices $i$ for which $M_i/M_{i-1}\cong \Delta^\fp_\lambda$ and the number of indices $i$ for which $M_i/M_{i-1}\cong \nb^\fp_\lambda$, respectively. 
The numbers $(M:\Delta^\fp_\lambda)$ and $(M:\nb^\fp_\lambda)$ are determined by~$\text{ch}M$, for $\la \in \Sigma_\fp^+.$

   We denote by $\mc O^{\Delta^\fp}$ (resp. $\mc O^{\nb^\fp}$) the full subcategory of modules that have a Verma flag (resp. a dual Verma flag).  Recall that the  ({\em indecomposable}) {\em tilting module} $T^\fp_\la$ studied in \cite[Theorem 3.5]{CCC} is an indecomposable module which lies in $\mc O^{\Delta^\fp}\cap \mc O^{\nb^\fp}$ with highest weight $\la -\rho$.

	
	\begin{lem}\cite[Corollary 3.8]{CCC}\label{[Corollary 3.8]} \label{lem::RiRe}   For any $\la, \mu \in \Sigma_\fp^+$, we have 
	\begin{align}
	&(T^\fp_\la: \Delta^\fp_{\mu}) = [\nabla^\fp_{-w_0^\fp\mu}: L_{-w_0^\fp\la}]
	\end{align} 
	\end{lem}

	We classify tilting modules whose characters are controlled by Kazhdan-Lusztig polynomials of type A Lie algebras as follows. 
	\begin{thm} \label{cor::ChP} \label{cor::chtypicalP} For any $\la,\mu \in \Sigma_\fp^+$, we have  \begin{align}
		&(T^\fp_\la :\nabla^\fp_{\mu}) = [\Delta^\fp_{-w_0^\fp\mu}: L_{-w_0^\fp\la}]. \label{Eq::tiltch}
		\end{align}
If $\la$ is $\mf p$-weakly-typical, then 
 \begin{align}
 &\ch T^\fp_\la = \sum_{\mu \in \h^\ast}[M^{\fp_\oa}_{-w_0^\fp\mu}: L^\oa_{-w_0^\fp\la}] \ch \nabla^\fp_{\mu}. \label{Eq::tiltch22}
 \end{align}  
 In particular, for any $\la \in \Sigma_\fp^+$ the following conditions are equivalent:
 \begin{itemize}
 	\item[(i).]$T^\fp_\la = \nabla^\fp_{\la}$.
 	\item[(ii).] $\la$  is $\mf p$-weakly-typical and $M^{\fp_\oa}_{-w_0^\fp\la}$ is projective in $\mc O(\mf g_\oa, \mf p_\oa)$. 
 \end{itemize}   
	\end{thm}
\begin{proof}
Recall from \cite[Proposition 3.4]{CCC} that there is a parabolic subalgebra $\hat \fp$ and a canonical contravariant equivalence $\bf D:\mc O(\mf g,  \fp)\rightarrow \mc O(\mf g, \hat \fp)$. By applying $\bf D$ together with \cite[Lemma 3.6]{CCC}  and  Lemma \ref{lem::RiRe} above, one deduces that 
	\begin{align}
	&(T^\fp_\la: \nb^\fp_{\mu})  =(T^{\hat \fp}_{-w_0\la}: \Delta^{\hat \fp}_{-w_0\mu}) = [\nb^{\hat \fp}_{w_0^{\hat \fp} w_0{\mu}}:L^{\hat \fp}_{w_0^{\hat \fp} w_0{\la}}]=[\Delta^{ \fp}_{-w_0 w_0^{\hat \fp} w_0{\mu}}:L^{\fp}_{-w_0 w_0^{\hat\fp} w_0{\la}}] =[\Delta^{ \fp}_{- w_0^{ \fp} {\mu}}:L^{\fp}_{- w_0^{\fp}{\la}}]
	\end{align} 
	This proves \eqref{Eq::tiltch} and also (ii) $\Rightarrow$ (i). Also, (\ref{Eq::tiltch22}) follows from (\ref{Eq::tiltch}) and Lemma \ref{prop::tyKac}.
	
	Now assume that $T_\la^\fp=\nb_\la^\fp$. Suppose on the contrary that $\la$ is not $\mf p$-weakly-typical. By \cite[Theorem 4.1]{CM}, there exists $\eta\in \h^\ast$ such that $L_{-w_0^\fp\la} ={\rm soc} K_\eta$. By Lemma \ref{lem::wt}, it follows that $\eta \neq -w_0^\fp\la$. Also,   there is an embedding  $\Lambda^{\rm top} \mf g_1 \otimes  L^\oa_{\eta} \hookrightarrow \Res L_{-w_0^\fp\la}$, which implies that $\mf p_\oa$ acts  locally finite  on $L^\oa_{\eta}$. In particular, we found that $\fp_\oa$ acts  locally finite   on $K_\eta$, namely, $\eta \in \Sigma_\fp^+$.  
	By \eqref{Eq::tiltch}, we have 
	$$(T^\fp_\la :\nabla^\fp_{-w_0^\fp\eta}) = [\Delta^\fp_{\eta}: L_{-w_0^\fp\la}] \geq  [K_{\eta}: L_{-w_0^\fp\la}]>0,$$ 
	which implies that $T_\la^\fp \neq \nb^\fp_\la$, a contradiction. This shows that $\la$ is $\fp$-weakly-typical. It follows from \eqref{Eq::tiltch22} and the BGG reciprocity that $M^{\fp_\oa}_{-w_0^\fp\la}$ is projective in $\mc O(\mf g_\oa, \mf p_\oa)$. This completes the proof. 	
\end{proof}

 \begin{rem}
 	Consider the special case $\mf p = \mf g_{\geq 0} :=\mf g_0\oplus \mf g_1$. Then for a given $\la \in \Sigma_{\mf g_{\geq 0}}^+$    we have  $$T_\la^{\mf g_{\geq 0}} =\nb_\la^{\g_{\geq 0}}\Leftrightarrow K_\la =L_\la \Leftrightarrow \la \text{ is typical.} $$ (cf. \cite{Ka2} and \cite[Lemma 3.4.1]{B+9}).
 \end{rem} 
 
We have the following descriptions for tilting characters in the full category $\mc O$.
 \begin{cor} Let $\la, \mu \in \h^\ast$. Then we have  
	 \begin{itemize}
	 		\item[(1)] $(T_\la :\nabla_{\mu}) = [\Delta_{-\mu}: L_{-\la}].$
		\item[(2)] $T_\la = \nabla_{\la}$ if and only if $\la$ is anti-dominant and $\mf b$-weakly-typical (i.e. anti-dominant and typical). 
		\item[(3)] 
		If $\la$ is integral and dominant, then $\ch T_\la = \sum_{\mu \in W \la} \ch \nabla_\mu.$
	\end{itemize}
\end{cor}

\begin{rem} It was shown in \cite[Section 8]{B+9} that indecomposable tilting modules in the category $\mc F_n$ of finite-dimensional modules are multiplicity-free. By Theorem \ref{cor::ChP} it is worthwhile to mention that for $\mf g=\pn$ with $n\geq 4$ there always exist $\la ,\mu\in \h^\ast$ such that $(T_\la:\Delta_\mu)>1$. 
\end{rem}

\subsection{Translation functors}
\subsubsection{Translation functors on $\mc O$} \label{Sect::transfun}  We recall the translation functor defined in \cite[Corollary 5.9]{CC}. Let $\C^{n|n}$ be the natural representation of $\gl(n|n)$. Then the exact functor $-\otimes \C^{n|n} :\mc O \rightarrow \mc O$  is decomposed as the direct sum of subfunctors $$-\otimes \C^{n|n} = \bigoplus_{a\in \C} (-\otimes \C^{n|n})_a,$$ according to eigenvalues $a$ of the {\em fake Casimir element} $\Omega$ (cf. \cite[Section 8.4]{Co}, \cite[Section 4.1]{B+9} and \cite[Section 2]{CP}), where $\theta_a:=(-\otimes \C^{n|n})_a : \cO \to\cO$ is the subfunctor of taking the $a$-eigenspace.

In this subsection we describe the action of $\theta_a$ on the Grothendieck groups of $\mc O^{\Delta}$ and $\mc O^\nb$. For weights $\la,\mu\in \h^\ast$ and $a\in \C$,  we  write $\la \dashrightarrow_a \mu$ if there exists some $1\leq i\leq n$ satisfying the following conditions: 
\begin{enumerate}
	\item $\mu_i=\la_i \pm 1$ and $\la_i =a$.
	\item $\mu_j=\la_j$, for any $j\neq i$. 
\end{enumerate} 
Also, we write $\la \rightsquigarrow_a \mu$ if there exists some $1\leq i\leq n$  satisfying the following conditions:  
\begin{enumerate}
	\item $\mu_i =\la_i+1$ and $\la_i =a$.
	\item $\mu_i =\la_i-1$ and $\la_i =a+2$.
	\item $\mu_j=\la_j$, for any $j\neq i$. 
\end{enumerate} 

\begin{lem}{\rm (}\cite[Proposition 5.7]{CC}{\rm)} \label{CCprop57}
	For any $a\in \mathbb C$, we have the following character formula: 
	\begin{align}
	&\emph{ch}\theta_a\Delta_\la = \sum_{\la \dashrightarrow_a \mu} \emph{ch}\Delta_\mu.\label{Eq::tra1} 
	\end{align}
\end{lem}
Also, we have the following rule:
\begin{lem} \label{CCprop58} For any $a\in \mathbb C$, we have the following character formula: 
	\begin{align}
	&\emph{ch}\theta_a\nabla_\la = \sum_{\la \rightsquigarrow_a \mu} \emph{ch}\nabla_\mu. \label{Eq::tra2}
	\end{align}
\end{lem}
\begin{proof}
	By \eqref{Eq::Delnb} we have $ \ch\nabla_\la= \sum_{\kappa \in S} \ch\Delta_{\la-\kappa}$. Since $\theta_a$ is an exact functor, by Lemma \ref{CCprop57} we have 
	\begin{align*}
	& \ch\theta_a\nabla_{\la} \\
	&=\sum_{j=1}^n \left(\sum_{\kappa_j=0,~\la_j=a}\ch \Delta_{\la-\kappa\pm \vare_j}+\sum_{\kappa_j=2,~\la_j=a+2}\ch \Delta_{\la-\kappa\pm \vare_j}\right)\\
	&=\sum_{j=1}^n \left(\sum_{\kappa_j=0,~\la_j=a}{\ch}\Delta_{\la+ \vare_j -\kappa} +\sum_{\kappa_j=2,~\la_j=a}{\ch}\Delta_{\la+ \vare_j-\kappa}
	+\sum_{\kappa_j=0,~\la_j=a+2}{\ch}\Delta_{\la- \vare_j -\kappa} +\sum_{\kappa_j=2,~\la_j=a+2}{\ch}\Delta_{\la- \vare_j -\kappa} \right) \\
	&=\sum_{\la \rightsquigarrow_a \mu} \ch\nabla_{\mu}.
	\end{align*} 
\end{proof}


\subsubsection{Translation functors on $\mc O^\fp$}\label{Sect::transfunOp}
As a consequence of Lemma \ref{CCprop57} and Lemma \ref{CCprop58}, we have analogous character formulae for parabolic categories. 

By definition, $\mc O^\fp$ is a full subcategory of $\mc O$. We denote the exact full embedding functor of categories $\mc O^\fp \hookrightarrow \mc O$. Its left adjoint functor   $Z^\fp:\mc O\rightarrow \mc O^\fp$ is the corresponding {\em Zuckerman functor}, taking the largest quotient inside $\mc O^\fp$. 
 In particular, $Z^\fp$ is a right exact functor.  We have $Z^\fp\Delta_\la \cong \Delta_\la^\fp$ for any $\la \in \Sigma^+_\fp$. We refer the reader to \cite[Section 9]{Hu08} for more details.

We recall the ordering defined in \ref{Eq::chord}. We have the following inequality of characters:
\begin{lem} \label{lem411}
	Let $M\in \mc O$ be given with the following  filtration, $$0=M_0\subseteq M_1\subseteq M_2 \subseteq \cdots \subseteq M_{\ell+1} =M.$$
	Then we have $\ch Z^\fp (M) \leq \sum\limits_{0\leq i\leq \ell} \ch Z^\fp(M_{i+1}/M_i).$ 
\end{lem}
\begin{proof}
	We prove this lemma by induction on $\ell$. Since $Z^\fp$ is right exact, we have exact sequence
	$$Z^\fp(M_\ell) \rightarrow Z^\fp(M) \rightarrow Z^\fp(M_{\ell+1}/M_\ell) \rightarrow 0.$$
	By induction hypothesis, we have 
	$$\ch Z^\fp(M)\leq \ch Z^\fp(M_\ell) +   \ch Z^\fp(M_{\ell +1}/M_\ell) \leq \sum_{0\leq i\leq  \ell-1}\ch Z^\fp(M_{i+1}/M_i) + \ch Z^\fp(M_{\ell +1}/M_\ell),$$ as desired. 
\end{proof}

\begin{prop}
Let $\la \in \Sigma_\fp^+$. For any $a\in \mathbb C$, we have the following character formulae:
	\begin{align}
	&\emph{ch}\theta_a\Delta^\fp_\la = \sum_{\substack{\mu \in \Sigma^+_\fp\\  \la \dashrightarrow_a \mu}} \emph{ch}\Delta^\fp_\mu.\label{Eq::tra3}  \\
	&\emph{ch}\theta_a\nabla^\fp_\la = \sum_{\substack{\mu \in \Sigma^+_\fp\\ \la \rightsquigarrow_a \mu }} \emph{ch}\nabla^\fp_\mu. \label{Eq::tra4}
	\end{align}    
\end{prop}
\begin{proof} By  \cite[Proposition 6.5]{Kn88}, we have  $\Delta_\la \otimes V\in \mc O^\Delta$ and $\Delta^\fp_\la \otimes V \in \mc O^{\Delta^\fp}$. 
Since $\theta_a$ is exact, the Zuckerman functor $Z^\fp$ gives rise to an epimorphism $\theta_a \Delta_\la \twoheadrightarrow \theta_a \Delta_\la^\fp \cong Z^\fp\theta_a \Delta_\la$. By Lemma \ref{lem411},  we have
\begin{align*}
(\theta_a \Delta_\la^\fp: \Delta_\mu^\fp) \leq (\theta_a \Delta_\la: \Delta_\mu),\end{align*} 
	 for any $\mu\in \Sigma_\fp^+$.
	 
Furthermore, by Pieri's rule we have 
\begin{align*}
\ch \Delta_\la \otimes V  =  & \quad \sum_{i=1}^n\ch \Delta_{\la +\vare_i} +\sum_{i=1}^n\ch \Delta_{\la -\vare_i},\\
\ch \Delta^\fp_\la \otimes V  =  & \sum_{\substack{1\leq i\leq n, \\ \la +\vare_i\in \Sigma^+_\fp}}\ch \Delta^\fp_{\la +\vare_i} +\sum_{\substack{1\leq i\leq n, \\ \la -\vare_i\in \Sigma^+_\fp}}\ch \Delta^\fp_{\la -\vare_i}.
\end{align*}
This means that  for any $\mu \in \Sigma^+_\fp$ we have 
\begin{align} \label{Eq::123}
&(\Delta^\fp_\la\otimes V: \Delta^\fp_\mu)=(\Delta_\la\otimes V: \Delta_\mu),
\end{align} which implies that $ (\theta_a\Delta^\fp_\la: \Delta^\fp_\mu)=(\theta_a \Delta_\la: \Delta_\mu)$ since $M=\bigoplus_{a\in \C}\theta_a M$ for any $M\in \mc O$. This completes the proof.

\end{proof} 
\vskip 0.5cm

\section{Character formulae of tilting modules of $\mc O^\fp$ for $\mf{pe}(3)$} \label{Sect::tranProj}

In this section we focus on the case when $\mf g=\mf{pe}(3)$.
 For simplicity, if  {$M\in \mc O^{\Delta^\fp}$}, we use the following expression in the remaining part of this article 
\begin{align*}
&M = \sum_{i \in I} k_i \Delta^\fp_{\la_i},
\end{align*}  instead of writing $\ch M = \sum_{i \in I} k_i \ch \Delta^\fp_{\la_i}$ for some coefficients $k_i\in \Z_{\geq 0}$. 
Similarly, we use analogous expression of characters in terms of characters of $\nabla^\fp$, Kac modules and irreducible modules, etc.

We will give the complete list of character formulae of tilting modules in $\mc O$.  
By $T^\fp_{a,b,c}$, we mean the tilting module $T^\fp_{a\vare_1 +b\vare_2+c\vare_3}$ with highest weight $a\vare_1 +b\vare_2+c\vare_3$.
Similarly, we define $\nabla^\fp_{a,b,c}$, $L_{a,b,c}$, etc. Also, we define $T_\la =T^{\mf b}_\la$, for any $\la \in \h^\ast$.

To simplify notation, we say a weight $\la  \in \h^\ast$ is  weakly typical  if $\la$ is $\mf b$-weakly-typical (i.e., $-\la$ is $\mf g_\oa$-weakly-typical),  which is equivalent to the condition $\prod_{\beta\in \Phi_\oa^+}(\langle\la,\beta\rangle   +1 ) \neq 0$.

\subsection{Character formulae of tilting modules for integral blocks in $\mc O$}

We first give the complete list of character formulae of tilting modules in $\mc O_\Z$. Note that the cases of weakly typical tilting modules have been given in Theorem \ref{cor::chtypicalP}. Moreover, observe that $T_{a,b,c}\otimes \mathbb C_{k\omega_3} \cong T_{a+k,b+k,c+k},$ for any $k\in\mathbb{C}$. 
Consequently, the character problem of an arbitrary tilting module in $\mc O_\Z$ is reduced to the following special cases: $T_{0,1,b}$, $T_{0,c,1}$ and $T_{a,0,1}$, where $a,b,c\in\Z$. In particular, the multiplicities $(T_\la:\Delta_\mu)$ are always at most one, for any $\la,\mu\in\h^\ast$.

 \subsubsection{Case $T_{0,1,b}$}
 \begin{thm} \label{thm::chTilt}
 	We have the following list of  character formulae of tilting modules.
 	\begin{align}
 	  T_{0,1,b} &=  \nabla_{0,1,b} + \nabla_{-1,0,b},~\text{for $b>2$}. \label{Eq::chT1}\\
 	  T_{0,1,1} &=  \nabla_{0,1,1} + \nabla_{-1,0,1}  +   \nabla_{-1,1,0}. \label{Eq::chT3} \\ 
 	  T_{0,1,b} &=\nabla_{0,1,b}+\nabla_{0,b,1}    +   \nabla_{b,0,1}  +   \nabla_{b,1,0} +  \nabla_{-1,0,b}   + \nabla_{-1,b,0}    +   \nabla_{b,-1,0}  +   \nabla_{b,0,-1}, ~\text{for $b<-1$}. \label{Eq::chT4} \\ 
 	   T_{0,1,-1} &= \nabla_{0,1,-1} +\nabla_{0,-1,1}   + \nabla_{-1,1,0} +\nabla_{-1,0,1} +  \nabla_{-1,0,-1} +  \nabla_{-1,-1,0}. \label{Eq::chT5} \\
 	    T_{0,1,0} &=\nabla_{0,1,0}   + \nabla_{0,0,1}    +   \nabla_{-1,0,0}.   \label{Eq::chT6} 
 	\end{align}
 \end{thm}
  
	\begin{rem}
		The character of $T_{0,1,2}\cong T_{-1,0,1}\otimes  \mathbb C_{\omega_3}$ follows from (\ref{Eq::chT996}). 
	\end{rem}

\subsubsection{Case $T_{0,c,1}$}
\begin{thm}\label{thm::chTilt2}
	We have the following list of  character formulae of   tilting modules.
	\begin{align}
	T_{0,c,1} &=  \nb_{0,c,1}+ \nb_{0,1,c} + \nb_{-1,c,0}  + \nb_{-1,0,c},~\text{for $c>1$}. \label{Eq::chT91}\\
	T_{0,1,1} &= \nb_{0,1,1} + \nb_{-1,0,1}  + \nb_{-1,1,0}. \label{Eq::chT92}\\
	T_{0,-1,1}	&=  \nb_{0,-1,1}  + \nb_{-1,0,1} + \nb_{-1,-1,0} . \label{Eq::chT93}\\
	T_{0,c,1} &=  \nb_{0,c,1}+ \nb_{c,0,1} + \nb_{-1,c,0}  + \nb_{c,-1,0},~\text{for $c<-1$}.\label{Eq::chT94}
	\end{align}
\end{thm}

\subsubsection{Case $T_{a,0,1}$}
\begin{thm}\label{thm::chTilt3}
	We have the following list of  character formulae of tilting modules.
	\begin{align}
	  T_{a,0,1} &=  \nb_{a,0,1}   +  \nb_{a,-1,0}  ,~\text{ for $a<-2$.} \label{Eq::chT991}\\
	  T_{a,0,1} &=  \nb_{a,0,1}  +   \nb_{0,a,1} +  \nb_{1,0,a}  +   \nb_{0,1,a}  +  \nb_{a,-1,0}  +  \nb_{-1,a,0} +  \nb_{0,-1,a}      +  \nb_{-1,0,a} ,~\text{ for $a>1$.} \label{Eq::chT992}\\
	  T_{1,0,1} &=  \nb_{1,0,1}   +  \nb_{0,1,1}  + \nb_{1,-1,0}   +  \nb_{0,-1,1} +  \nb_{-1,0,1}   +   \nb_{-1,1,0}  . \label{Eq::chT993}\\
	  T_{0,0,1} &=  \nb_{0,0,1}   +  \nb_{0,-1,0}   +  \nb_{-1,0,0}.  	\label{Eq::chT994} \\
	  T_{-2,0,1} &=  \nb_{-2,0,1}   +  \nb_{-2,-1,0}   +  \nb_{-3,-2,0}.  	\label{Eq::chT995} \\
	  T_{-1,0,1} &=  \nb_{-1,0,1}   +  \nb_{-1,-1,0}    +   \nb_{-2,-1,1}+  \nb_{-3,-1,0}   +  \nb_{-2,-1,-1}   +  \nb_{-3,-2,-1}  	\label{Eq::chT996} 
	\end{align}
\end{thm}

The proofs of these theorems are based on straightforward computations and case-by-case discussions.   The details can be found in the appendix.

\subsection{Character formulae of tilting modules for non-integral blocks in $\mc O$}
In this subsection, we consider tilting modules of non-integral weights. Let $\la =\la_1\vare_1 +\la_2\vare_2 + \la_3 \vare_3 \in \h^\ast$ and $1\leq i\leq 2$ be such that   $\langle\la, \alpha_i\rangle \notin \Z$.  Then it follows from \cite[Proposition 4.4]{CC} that $\mc O_\la \cong  \mc O_{s_{\alpha_i}\la}$ (see also \cite[Proposition 8.6]{CoM1}).  We may assume that   $\langle \la,\alpha_1\rangle \in \Z$ and $\langle \la,\alpha_2\rangle \notin \Z$. 
As a result, the tilting character problem reduces to the cases in the following theorem. 
\begin{thm} \label{thm::nonintegralchp3}
	We have the following list of  character formulae of tilting modules for $b\in \Z$ and $c\notin \Z$: 
	\begin{align*}
	&T_{0,b,c} = \nb_{0,b,c}, ~\text{ for $b\geq 2$.}\\
	& T_{0,1,c} = \nb_{0,1,c}+\nb_{-1,0,c}. \\
	& T_{0,0,c} = \nb_{0,0,c}.\\
	&T_{0,b,c} = \nb_{0,b,c}+\nb_{b,0,c}, ~\text{ for $b\leq -1$.}\\
	\end{align*}
\end{thm}
\begin{proof}
	By Theorem \ref{cor::ChP} it remains to show that $ T_{0,1,c} = \nb_{0,1,c}+\nb_{-1,0,c}$. Observe that $T_{-1,1,c} =\nb_{-1,1,c}$ and thus $\theta_{-1}T_{-1,1,c} = \nb_{0,1,c}+\nb_{-1,0,c} = T_{0,1,c}.$
\end{proof} 

 Suppose that $\la = \la_1\vare_1+\la_2\vare_2+c\vare_3$ and $\mu=\mu_1\vare_1+\mu_2\vare_2+c\vare_3$ are $\mf{pe}(3)$-weights for some $\la_1,\la_2,\mu_1,\mu_2\in \Z, c\notin \Z$. Define the $\mf{pe}(2)$-weights $\la' :=\la_1\vare_1+\la_2\vare_2$ and $\mu' :=\mu_1\vare_1+\mu_2\vare_2$. The following corollary is a direct consequence of Theorem \ref{thm::nonintegralchp3}. 
\begin{cor} Let $\mc O_n$ denotes the cateogry $\mc O$ for $\mf{pe}(n)$. Let $P_{\la}$ and $P_{\la'}$ denote the projective covers of $L_{\la}$ and $L_{\la'}$ in $\mc O_3$ and $\mc O_2$, respectively. Then we have 
\[{\rm dim \, Hom}_{\mc O_3}(P_\la, P_\mu) = {\rm dim \, Hom}_{\mc O_2}(P_{\la'}, P_{\mu'}).\]   
\end{cor}

\begin{rem}
	  In \cite{CMW} Cheng, Mazorchuk and Wang proved   that non-integral blocks of
	category $\mc O$ for general linear Lie superalgebras are equivalent to integral blocks in $\mc O$ for direct sums of general linear Lie superalgebras with strictly smaller ranks. 
	A similar phenomenon was observed in
	\cite[proof of Theorem 6.5]{CC}. That is,  non-integral $\mf{pe}(2)$-blocks are equivalent to integral blocks for $\mf{pe}(1)\oplus \mf{pe}(1)$. 
	However, as pointed out in \cite[Remark 6.6]{CC}, in the case of $\mf{pe}(2)$, there is no parabolic
	induction as in \cite{CMW} which directly provides such equivalence.
	
	Observe that Theorem \ref{thm::nonintegralchp3} and its proof imply the isomorphism between the Grothendieck group for the non-integral $\mf{pe}(3)$-block $\mc O_\la$ and that for the integral blocks for $\mf{pe}(2)\oplus \mf{pe}(1).$ We should mention that a similar argument in the proof of Theorem \ref{thm::nonintegralchp3} can be generalized to an arbitrary $\mf{pe}(n)$. For instance, let $\la\in \h^\ast$ be a weight of arbitrary  $\mf{pe}(n)$. If $\langle\la, \alpha_1\rangle \in \Z$ and $\langle\la, \vare_{i}-\vare_j\rangle   \notin \Z$, for any $3\leq i<j\leq n$, 
	then the Grothendieck group for the non-integral $\mf{pe}(n)$-block $\mc O_\la$ and that for integral blocks for $\mf{pe}(2)\oplus \mf{pe}(1)^{\oplus n-2}$ are naturally isomorphic. 
\end{rem}

\subsection{Character formulae of tilting modules for other parabolic categories}
Recall that tilting modules coincide with projective modules in the category of finite-dimensional modules for $\pn$. By    Section \ref{Sect::tranProj} and \cite[Section 3]{Ch} (see also \cite{B+9}), we are left with the problem of tilting characters in   $\mc O^\fp$ for $\fp_\oa \neq \g_\oa,~\h$.

For any parabolic subalgebra $\mf s\subseteq \mf g$, let $\Phi^+(\mf s)$ denote the set of positive roots of $\mf s$. We now set $\fp$ to be the  parabolic subalgebra determined by  $\Phi^+(\mf p_\oa) = \Phi^+(\mf b_\oa)\cup \{ -\alpha_1 \}$ and  $\mf p_\ob =\mf g_{1}$.
The  corresponding Levi subalgebra $\mf l\cong \gl(2)\oplus \gl(1)$  have roots $\{\pm\alpha_1\}$.  

Recall the canonical contravariant equivalence $\bf D:\mc O(\mf g, \mf p)\rightarrow \mc O(\mf g, \hat \fp)$ in \cite[Proposition 3.4]{CCC}, where $\hat \fp$ is a  parabolic subalgebra determined by  $\Phi^+(\mf p_\oa) = \Phi^+(\mf b_\oa)\cup \{ -\alpha_2 \}$ and  $\mf p_\ob =\mf g_{-1}$. The corresponding Levi subalgebra $ \hat{\mf l}$ is generated by root vectors of even root $\{ \pm\alpha_2 \}$.
Using the algorithm in \cite[Section 2.2]{PS89},
the problem of general tilting characters in $\mc O^\fp$ is reduced to the following special cases.

\begin{thm}\label{tiltcom} Let $\mf p$ be the parabolic subalgebra given as above. Then  
	we have the following list of  character formulae of tilting modules of non-$\mf p$-weakly-typical highest weights: 
	\begin{align*}
		&T^\fp_{1,0,a} =\nb^\fp_{1,0,a} + \nb^\fp_{0,-1,a}, ~\text{for $a\geq 3$.} \\
		&T^\fp_{1,0,2}= \nb^\fp_{1,0,2}+\nb^\fp_{0,-1,2}+\nb^\fp_{0,-2,1}+\nb^\fp_{0,-1,0}+\nb^\fp_{-1,-2,0}.\\
	&T^\fp_{1,0,1} = \nb^\fp_{1,0,1} +\nb^\fp_{0,-1,1}+\nb^\fp_{1,-1,0}.\\
	&T^\fp_{1,0,0} = \nb^\fp_{1,0,0}+\nb^\fp_{0,-1,0}.\\
	&T^\fp_{1,0,-1} = \nb^\fp_{1,0,-1}+\nb^\fp_{1,-1,0}+\nb^\fp_{0,-1,-1}.\\
	&T^\fp_{1,0,a} = \nb^\fp_{1,0,a}+\nb^\fp_{1,a,0}+\nb^\fp_{0,-1,a}+\nb^\fp_{0,a,-1},~\text{for $a\leq -2$.}\\
	&T^\fp_{1,a,2} = \nb^\fp_{1,a,2}+\nb^\fp_{0,a,1}, ~\text{for $a<-1$.}\\
	&T^\fp_{1,-1,2}=\nb^\fp_{1,-1,2}+\nb^\fp_{-1,-2,1}+\nb^\fp_{0,-1,1}.
	\end{align*}
\end{thm}

 \begin{rem}
 	Similar to the case of the full category $\mc O$, by tensoring with the one-dimensional representation $\C_{k\omega_3}$, the characters provided in Theorem~\ref{tiltcom} indeed give a complete list of tilting characters in $\mc O_\Z^\fp$.
 \end{rem}

    \section{Appendix} \label{Section::app}
    In the appendix, a weight $\la=\la_1\ep_1+\la_2\ep_2+\la_3\ep_3\in \mf h^*$ will be written as $(\la_1,\la_2,\la_3)$.
    Moreover, the following standard facts will be frequently used:    
    \begin{enumerate}
    	\item[(1)] If  $\mu \uparrow w \mu $  and $[K_{\mu}:L_{\la}]>0$, then  $$(T_{-\la}: \nabla_{ -w\mu})= [\Delta_{w \mu}:L_{\la}] =(\Delta_{w\mu}:K_{\mu})[K_{\mu}:L_{\la}]>0.$$ 
    	\item[(2)] If  $L_{\la}=L^{\mf b^r}_{\mu -2\omega_3}$, then  $$(T_{-\la}: \nabla_{ -w\mu})= [\Delta_{w \mu}:L_{\la}]>0  \quad \text{for any}\quad \mu \uparrow w \mu. $$ 
    \end{enumerate}

    \subsection{Proof of Theorem \ref{thm::chTilt}}
    The proof is based on case-by-case discussions. \vskip 0.1cm
    {{\bf Case (I): proof of} \eqref{Eq::chT1}}.\\ Observe that $T_{-1,1,b} = \nb_{-1,1,b}$ since $(-1,1,b)$ is weakly typical and anti-dominant for $b>2$. 
    Applying the functor $\theta_{-1}$, we obtain $\theta_{-1}T_{-1,1,b} = \nb_{0,1,b}+\nb_{-1,0,b}$. 
    By Proposition \ref{cor::chtypicalP}, $T_{0,1,b}$ is a direct summand of $\theta_{-1}T_{-1,1,b}$. 
    It follows  from Corollary \ref{cor::oddlinkage} that $$(T_{0,1,b}: \nb_{-1,0,b}) = [\Delta_{1,0,-b} : L_{0,-1,-b}] > 0.$$

    {{\bf Case (II): proof of} \eqref{Eq::chT3}}.\\
    Observe that $(-1,1,1)$ is weakly typical and anti-dominant. Applying $\theta_{-1}$, we have $$\theta_{-1} T_{-1,1,1} = \nb_{0,1,1}+\nb_{-1,0,1}+\nb_{-1,1,0}.$$ 
    Note that $T_{0,1,1}$ is a direct summand of $\theta_{-1}T_{-1,1,1}$. By Corollary \ref{cor::oddlinkage}, we have $(T_{0,1,1}: \nb_{-1,0,1})>0$. 
    Finally, $T_{-1,1,0}\neq \nb_{-1,1,0}$ since $(-1,1,0)$ is not anti-dominant. We conclude that $\theta_{-1} T_{-1,1,1} = T_{0,1,1}$.

    {{\bf Case (III): proof of} \eqref{Eq::chT4}}.\\ 
    Since $(-1,1,b)$ is weakly typical for $b<-1$, we have $T_{-1,1,b} = \nb_{-1,1,b} +\nb_{-1,b,1}+\nb_{b,-1,1} +\nb_{b,1,-1}$. 
    Applying $\theta_{-1}$, we have 
    \begin{align}
    &\theta_{-1} T_{-1,1,b} = \nb_{0,1,b} +\nb_{-1,0,b} +\nb_{0,b,1} +\nb_{-1,b,0} +\nb_{b,0,1} +\nb_{b,-1,0}+\nb_{b,0,-1} +\nb_{b,1,0}. 
    \end{align} 
    
    Observe that $T_{0,1,b}$ is a direct summand of $\theta_{-1} T_{-1,1,b}$. By Corollary \ref{cor::oddlinkage} and Proposition \ref{prop::Jan}, it follows that  
    $(T_{0,1,b} : \nb_\mu) = [\Delta_{-\mu} : L_{0,-1,-b}] >0$  for any $\mu$ such that $(1,0,-b) \uparrow -\mu$. Consequently, we have $T_{0,1,b} = \theta_{-1}T_{-1,1,b}.$

    {{\bf Case (IV): proof of} \eqref{Eq::chT5}}.\\	  Observe that $(-1,1-1)$ is weakly typical. Therefore 
    \begin{align}
    &T_{-1,1,-1} = \nb_{-1,1,-1} +\nb_{-1,-1,1}.
    \end{align} Applying $\theta_{-1}$, we have 
    \begin{align}
    &\theta_{-1} T_{-1,1,-1} = \nb_{0,1,-1} +\nb_{-1,0,-1} +\nb_{-1,1,0}  +\nb_{0,-1,1} +\nb_{-1, 0 ,1} +\nb_{-1,-1,0}.
    \end{align}
    
    Observe that $T_{0,1,-1}$ is a direct summand of $\theta_{-1} T_{-1,1,-1}$. As a consequence of Corollary \ref{cor::oddlinkage}, we have $(T_{0,1,-1}: \nb_{-1,0,-1}) >0$. By Proposition \ref{prop::Jan}, we may conclude that  $T_{0,1,-1} =\theta_{-1} T_{-1,1,-1}.$

    {{\bf Case (V): proof of} \eqref{Eq::chT6}}.\\  We firstly solve for the $\nabla$-flag in $T_{-1,1,0}$. Observe that $T_{-2,1,0} = \nb_{-2,1,0} +\nb_{-2,0,1}$ since $(-2,1,0)$ is weakly typical. Applying $\theta_{-2}$, we have 
    \begin{align}\label{6.4}
    &\theta_{-2} T_{-2,1,0} = \nb_{-1,1,0} +\nb_{-2,1,-1} +\nb_{-1,0,1} +\nb_{-2,-1,1}.
    \end{align} 
    Using odd reflections, we have $L_{1,-1,0} =\text{soc}K_{2,-1,1}$, which implies that $(T_{-1,1,0}: \nb_{-2,1,-1}) >0.$ Moreover, since $(-2,-1,1)$ is not weakly typical, we have $T_{-2,-1,1} \neq \nb_{-2,-1,1}$. We may conclude that $\theta_{-2}T_{-2,1,0} = T_{-1,1,0}$. 
    
    Applying the functor $\theta_{-1}$ to (\ref{6.4}), we have
    \begin{align}
    &\theta_{-1} T_{-1,1,0} = \nb_{0,1,0} +\nb_{-1,0,0} +\nb_{-2,0,-1} +\nb_{-2,1,0} +\nb_{0,0,1} +\nb_{-1,0,0} +\nb_{-2,0,1} +\nb_{-2,-1,0}. \label{Eq::Eq1}
    \end{align} Note that $T_{0,1,0}$ is a direct summand of $\theta_{-1} T_{-1,1,0}$. 
    
    On the other hand, by applying the functor $\theta_{-1}$ to \eqref{Eq::chT5}, we have 
    \begin{align}
    &\theta_{-1} T_{0,1,-1} = 2(\nb_{0,1,0} +\nb_{0,0,-1} +2\nb_{-1,0,0} +\nb_{0,0,1} + \nb_{0,-1,0}). \label{Eq::Eq2}
    \end{align}  
    Observe that $2T_{0,1,0}$ is a direct summand of $\theta_{-1} T_{0,1,-1} $. Combining \eqref{Eq::Eq1} and \eqref{Eq::Eq2}, we may conclude that $$T_{0,1,0} = \nb_{0,1,0} +a\nb_{-1,0,0}+b\nb_{0,0,1},$$  for some $a\in\{1,2\}$ and $b\in\{0,1\}$ by Corollary \ref{cor::oddlinkage}. Since $(T_{0,1,0}: \nb_{0,0,1}) = [\Delta_{0,0,-1}: L_{0,-1,0}] >0$, we must have $b=1$.  
    Suppose on that $a=2$. It follows from \eqref{Eq::Eq2} that $\nb_{0,0,-1}  +\nb_{0,-1,0}$ is the character of either $T_{0,-1,0}$ or $T_{0,-1,0}\oplus T_{0,0,-1}$, namely, it leads to $T_{0,-1,0} = \nb_{0,-1,0} +\nb_{0,0,-1}$, which   contradicts to Corollary \ref{cor::oddlinkage}. This completes the proof.

    \subsection{Proof of Theorem \ref{thm::chTilt2}}
    The proof is based on case-by-case discussions. \vskip 0.1cm
    {{\bf Case (I): proof of} \eqref{Eq::chT91}}.\\ 
    For $c>1$, $(-1,c,1)$ is weakly typical and hence $T_{-1,c,1} = \nb_{-1,c,1} +\nb_{-1,1,c}$. Consequently we have 
    \begin{align}
    & \theta_{-1} T_{-1,c,1} = \nb_{0,c,1}+\nb_{-1,c,0} +\nb_{0,1,c} +\nb_{-1,0,c}.
    \end{align}  
    Observe that $T_{0,c,1}$ is a direct summand of  $\theta_{-1} T_{-1,c,1}$. 
    Note that 
    $$(T_{0,c,1}: \nb_{0,1,c}) = [{ \Delta_{0,-1,-c}}: L_{0,-c,-1}]>0.$$
    Moreover,	we have $$(T_{0,c,1}: \nb_{-1,c,0}) = [{ \Delta}_{1,-c,0} : L_{0,-c,-1}] = ({ \Delta}_{1,-c,0} : K_{1,-c, 0})[K_{1,-c,0} : L_{0,-c,-1}].$$ 
    Using odd reflections, we find that 
    $$L_{0,-c,-1}  =  {L_{-2,-c-1,-1} = L^{\mf b^r}_{3, -c-3,-2}=}\text{soc} K_{1,-c,0},$$ which implies that $(T_{0,c,1}: \nb_{-1,c,0})>0$. 
    Finally, $T_{-1,0,c}\neq \nb_{-1,0,c}$ since $(-1,0,c)$ is not weakly typical. 
    This shows that $\theta_{-1} T_{-1,c,1} = T_{0,c,1}$, as desired.

    {{\bf Case (II): proof of} \eqref{Eq::chT92}}.\\  
    Observe that $(-1,1,1)$ is weakly typical and anti-dominant. Therefore $T_{-1,1,1} =\nb_{-1,1,1}$, and hence 
    \begin{align}
    &\theta_{-1} \nb_{-1,1,1} = \nb_{0,1,1} +\nb_{-1,0,1} +\nb_{-1,1,0}. 
    \end{align} 
    Note that $T_{0,1,1}$ is a direct summand of $\theta_{-1}T_{-1,1,1}$. 
    By Corollary \ref{cor::oddlinkage}, we have $(T_{0,1,1}: \nb_{-1,0,1}) >0.$ 
    Finally, since $(-1,1,0)$ is not anti-dominant, we have $T_{-1,1,0} \neq \nb_{-1,1,0}$, as desired.

    {{\bf Case (III): proof of} \eqref{Eq::chT93}}.\\   
    Observe that $(-1,-1,1)$ is weakly typical and anti-dominant. Hence we have $T_{-1,-1,1} = \nb_{-1,-1,1}$ and
    \begin{align}
    &\theta_{-1}T_{-1,-1,1}  =\nb_{0,-1,1} +\nb_{-1,0,1}+\nb_{-1,-1,0}.	
    \end{align}  Note that $T_{0,-1,1}$ is a direct summand of $\theta_{-1}T_{-1,-1,1}$. Also, we have  $(T_{0,-1,1} : \nb_{-1,0,1}) > 0 $. The proof now follows from the fact that  $T_{-1,-1,0} \neq \nb_{-1,-1,0}$ since $(-1,-1,0)$ is not weakly typical.

    {{\bf Case (IV): proof of} \eqref{Eq::chT94}}.\\  
    For $c<-1$, $(-1,c,1)$ is weakly typical. We claim that $T_{0,c,1} = \theta_{-1} T_{-1,c,1}$. Observe that  
    \begin{align}
    &\theta_{-1} T_{-1,c,1} = \theta_{-1}(\nb_{-1,c,1} +\nb_{c,-1,1}) = \nb_{0,c,1}+\nb_{-1,c,0}+\nb_{c,0,1}+\nb_{c,-1,0}.
    \end{align} Therefore $T_{0,c,1}$ is a direct summand of $\theta_{-1} T_{-1,c,1}$. 
     {Using odd reflections, we find that }
    $$L_{0,-c,-1} = L^{\mf b^r}_{-1,-c-2,-2} =\text{soc}K_{1,-c,0},$$  
    which implies that $(T_{0,c,1}: \nb_{-1,c,0})>0$. 
    Finally, $T_{c,-1,0}\neq \nb_{c,-1,0}$ since $(c,-1,0)$ is not weakly typical. This completes the proof.

    \subsection{Proof of Theorem \ref{thm::chTilt3}}
     The proof is based on case-by-case discussions. \vskip 0.1cm
    {{\bf Case (I): proof of} \eqref{Eq::chT991}}.\\   
     {For $a<-2$, $(a,-1,1)$ is weakly typical and hence $T_{a,-1,1} = \nb_{a,-1,1}$. 
    	Applying $\theta_{-1}$, we have}
    \begin{align}
    &\theta_{-1} T_{a,-1,1} = \nb_{a,0,1}   +  \nb_{a,-1,0}.
    \end{align} 
    Note that $ T_{a,0,1} $ is a direct summand of $\theta_{-1} T_{a,-1,1}$. 
    By Corollary \ref{cor::oddlinkage}, the statement follows.

    {{\bf Case (II): proof of} \eqref{Eq::chT992}}.\\   
    For $a>1$, $(a,-2,1)$ is weakly typical and we have 
    \begin{align}
    &T_{a,-1,1} = \nb_{a,-1,1} +\nb_{-1,a,1} +\nb_{-1,1,a} +\nb_{1,-1,a}.
    \end{align}
    Applying $\theta_{-1}$, we have 
    \begin{align}
    &\theta_{-1}T_{a,-1,1} = \nb_{a,0,1} +\nb_{a,-1,0} +\nb_{0,a,1} +\nb_{-1,a,0} +\nb_{0,1,a} +\nb_{-1, 0, a} +\nb_{0,-1,a} +\nb_{1,0,a}.
    \end{align} 
    Note that $T_{a,0,1}$ is a direct summand of $\theta_{-1}T_{a,-1,1}$. By Corollary \ref{cor::oddlinkage} and Proposition \ref{prop::Jan}, we have $$(T_{a,0,1}: \nb_{\mu}) = [\Delta_{-\mu}: L_{-a,0,-1}] >0, $$  {for any $\mu$ such that $(-a,1,0) \uparrow  -\mu$}, and the statement follows.

    {{\bf Case (III): proof of} \eqref{Eq::chT993}}.\\     
    We start with $T_{1,-1,1} = \nb_{1,-1,1} +\nb_{-1,1,1}$. Applying $\theta_{-1}$, we have 
    \begin{align}
    &\theta_{-1} T_{1,-1,1} = \nb_{0,-1,1} +\nb_{1,0,1} +\nb_{1,-1,0} +\nb_{0,1,1} +\nb_{-1,0,1} +\nb_{-1,1,0}.
    \end{align} 
    Note that $T_{1,0,1}$ is a direct summand of $\theta_{-1} T_{1,-1,1}$. 
    Using the odd reflections, we find that 
    \begin{align}
    &L_{-1,0,-1} = L^{\mf b^r}_{-3,-1,-2} =\text{soc}K_{-1,1,0}.
    \end{align} Therefore $$(T_{1,0,1}: \nb_{\mu}) = [\Delta_{-\mu}: L_{-1,0,-1}] >0, $$ for any $\mu$ such that $(-1,1,0) \uparrow  -\mu$, and the statement follows.

    {{\bf Case (IV): proof of} \eqref{Eq::chT994}}.\\     
    Applying $\theta_{-1}$  {to (\ref{Eq::chT93})}, we have 
    \begin{align}
    &\theta_{-1} T_{0,-1,1} = 2(\nb_{0,0,1} +\nb_{0,-1,0} +\nb_{-1,0,0}).
    \end{align}
    Note that $2T_{0,0,1}$  is a direct summand of $\theta_{-1} T_{0,-1,1}$. By Corollary \ref{cor::oddlinkage}, we have $(T_{0,0,1}: \nb_{0,-1,0}) >0$. Since $(-1,0,0)$ is not weakly typical, we have $T_{-1,0,0} \neq \nb_{-1,0,0}$, as desired.

    {{\bf Case (V): proof of} \eqref{Eq::chT995}}.\\    
    By \eqref{Eq::chT991}, we have $T_{-3,0,1} = \nb_{-3,0,1} +\nb_{-3,-1,0}$.  Applying the functor $\theta_{-3}$, we have  
    \begin{align}
    &\theta_{-3} T_{-3,0,1} =  \nb_{-2,0,1}   +  \nb_{-2,-1,0}   +  \nb_{-3,-2,0}.
    \end{align}Note that $T_{-2,0,1}$ is a direct summand of $\theta_{-3} T_{-3,0,1}.$  By Corollary \ref{cor::oddlinkage}, we have $(T_{-2,0,1}: \nb_{-2,-1,0})>0.$ 
    Since $(-3,-2,0)$ is not weakly typical, we have $T_{-3,-2,0} \neq \nb_{-3,-2,0}$, and the statement follows.

    {{\bf Case (VI): proof of} \eqref{Eq::chT996}}.\\   Since $(-3,-1,1)$ is weakly typical and anti-dominant, we may conclude that $T_{-3,-1,1} = \nb_{-3,-1,1}$. 
     Applying  $\theta_{-3}$, we have 
    \begin{align}\label{6.18}
    &\theta_{-3} T_{-3,-1,1} =\nb_{-2,-1,1} +\nb_{-3,-2,1}.
    \end{align} 
    Note that $T_{-2,-1,1}= \theta_{-3} T_{-3,-1,1}$ by Corollary \ref{cor::oddlinkage}. 
    Applying $\theta_{-1}$ to (\ref{6.18}), we have 
    \begin{align}
    &\theta_{-1} T_{-2,-1,1} = \nb_{-2,0,1} +\nb_{-2,-1,0} +\nb_{-3,-2,0},
    \end{align} which implies that  $\theta_{-1} T_{-2,-1,1}=T_{-2,0,1}$ by Corollary \ref{cor::oddlinkage} and the fact that $T_{-3,-2,0} \neq \nb_{-3,-2,0}$. 
     By Corollary \ref{cor::oddlinkage} again, we may conclude that 
    \begin{align}
    &T_{-1,0,1} = \nb_{-1,0,1}   +  \nb_{-1,-1,0}   +   \nb_{-2,-1,1}  +  a\nb_{-2,-1,-1}  +  b\nb_{-3,-1,0}    +  c\nb_{-3,-2,-1},  
    \end{align} for some $a,b,c \in \{0,1\}$.
    Now we have $a=1, b=1$ by \eqref{Eq::chT92} and \eqref{Eq::chT995}, respectively. Finally, note that $T_{-3,-2,-1} \neq \nb_{-3,-2,-1}$. This completes the proof.



\end{document}